\documentclass[twoside,11pt]{article}

\usepackage{jmlr2e}
\usepackage{mathtools}
\usepackage{lscape}
\usepackage{subcaption}
\usepackage{suffix}
\usepackage{color}
\usepackage{enumitem}
\usepackage{hyperref}

\usepackage{color,hyperref}
\definecolor{darkblue}{rgb}{0.0,0.0,0.5}
\definecolor{black}{rgb}{0.0,0.0,0.0}
\hypersetup{colorlinks,breaklinks,
            linkcolor=darkblue,urlcolor=darkblue,
            anchorcolor=darkblue,citecolor=darkblue}

\allowdisplaybreaks

\newcommand{\E}{{\bf E}}

\newcommand{\Nor}{{\cal N}}  

\newcommand{\ben}{\begin{enumerate}}
\newcommand{\een}{\end{enumerate}}
\newcommand{\beq}{\begin{equation}}
\newcommand{\eeq}{\end{equation}}

\newcommand{\argmin}{\operatornamewithlimits{argmin}}

\DeclarePairedDelimiterX\MeijerM[3]{\lparen}{\rparen}%
{\begin{smallmatrix}#1 \\ #2\end{smallmatrix}\delimsize\vert\,#3}

\newcommand\MeijerG[8][]{%
  G^{\,#2,#3}_{#4,#5}\MeijerM[#1]{#6}{#7}{#8}}

\WithSuffix\newcommand\MeijerG*[7]{G^{\,#1,#2}_{#3,#4}\MeijerM*{#5}{#6}{#7}}

\DeclarePairedDelimiterX\pFqM[3]{\lparen}{\rparen}%
{\begin{smallmatrix}#1 \\ #2\end{smallmatrix}\delimsize\vert\,#3}

\newcommand\pFq[6][]{%
  {}_{#2}F_{#3}\pFqM[#1]{#4}{#5}{#6}}

\WithSuffix\newcommand\pFq*[5]{{}_{#1}F_{#2}\pFqM*{#3}{#4}{#5}}

\newtheorem{theorem}{Theorem}
\numberwithin{theorem}{section}

\numberwithin{Def}{section}

\numberwithin{remark}{section}

\newtheorem{proposition}{Proposition}
\numberwithin{proposition}{section}
\newtheorem{lemma}{Lemma}
\numberwithin{lemma}{section}
\newtheorem{Cor}{Corollary}
\numberwithin{Cor}{section}

\newtheorem{result}{Result}

\usepackage{booktabs,array}

\newcount\rowc

\makeatletter
\def\ttabular{%
\hbox\bgroup
\let\\\cr
\def\rulea{\ifnum\rowc=\@ne \hrule height 1.0pt \fi}
\def\ruleb{
\ifnum\rowc=1\hrule height 1.0pt  \else
\ifnum\rowc= 3  \hrule height 0.5pt \else
\ifnum\rowc= 5  \hrule height 0.5pt \else
\ifnum\rowc= 7  \hrule height 0.5pt \else
\ifnum\rowc= 9  \hrule height 0.5pt \else
\ifnum\rowc= 11  \hrule height 0.5pt 
  \else \hrule height 0pt
\fi\fi\fi\fi\fi\fi}
\valign\bgroup
\global\rowc\@ne
\rulea
\hbox to 7em{\strut \hfill##\hfill}%
\ruleb
&&%
\global\advance\rowc\@ne
\hbox to 7em{\strut\hfill##\hfill}%
\ruleb
\cr}
\def\endttabular{%
\crcr\egroup\egroup}

\graphicspath{{../../figures/}{../figures/}{./figures/}{./}}
\renewcommand{\E}{\mathbb E}

\setcitestyle{citesep={;}}

\graphicspath{{../../figures/}{../figures/}{./figures/}{./}}

\jmlrheading{X}{XXXX}{XX--XX}{XX/XX}{XX/XX}{bhadraetal}{Anindya Bhadra, Jyotishka Datta, Yunfan Li, Nicholas G. Polson, and Brandon Willard}

\ShortHeadings{Prediction Risk}{Bhadra et al.}
\firstpageno{1}

\begin{document}

\title{Prediction Risk for the Horseshoe Regression}

\author{\name Anindya Bhadra \email bhadra@purdue.edu \\
       \addr Department of Statistics\\
       Purdue University\\
       West Lafayette, IN 47907, USA
       \AND
       \name Jyotishka Datta \email jd033@uark.edu\\
       \addr Department of Mathematical Sciences\\
       University of Arkansas\\
       Fayetteville, AR 72701, USA
       \AND
       \name Yunfan Li \email li896@purdue.edu\\
       \addr  Department of Statistics\\
       Purdue University\\
       West Lafayette, IN 47907, USA
        \AND
       \name Nicholas G. Polson \email ngp@chicagobooth.edu\\
       \addr Booth School of Business\\
       University of Chicago\\
       Chicago, IL 60637, USA
        \AND
       \name Brandon Willard \email brandonwillard@gmail.com\\
        \addr Booth School of Business\\
       University of Chicago\\
       Chicago, IL 60637, USA}

\editor{}

\maketitle

\begin{abstract}
We show that prediction performance for global-local shrinkage regression can overcome two major difficulties of global shrinkage regression:~(i) the amount of relative shrinkage is monotone in the singular values of the design matrix and (ii) the shrinkage is determined by a single tuning parameter. Specifically, we show that the horseshoe regression, with heavy-tailed component-specific local shrinkage parameters, in conjunction with a global parameter providing shrinkage towards zero, alleviates both these difficulties and consequently, results in an improved risk for prediction. Numerical demonstrations of improved prediction over competing approaches in simulations and in a pharmacogenomics data set confirm our theoretical findings.
\end{abstract}


\newpage
\section{Introduction}\label{sec:intro}
We develop theoretical results on prediction risk in the high-dimensional linear regression model
\begin{eqnarray}
y=X\beta + \epsilon,\label{eq:reg1}
\end{eqnarray}
where $y\in \mathbb{R}^n, X \in \mathbb{R}^{n\times p}, \beta \in \mathbb{R}^p$, $\epsilon \sim \Nor(0, \sigma^2\mathbf{I}_n)$ with $p>n$, and the design matrix $X$ is assumed fixed. Let $\hat\beta$ denote the estimate of $\beta$ based on the observed data $y$ and design $X$. Let $y^{*}$ denote a future observation generated from the same model, independent of $y$. Define the quadratic predictive risk
\begin{eqnarray}
R = \E_{y^{*}, y\mid  X, \beta} (y^* - X\hat\beta)^2, \label{eq:riskdef}
\end{eqnarray}
where the subscript denotes that the expectation is with respect to the data generating distribution, with $X$ and $\beta$ held fixed. We focus on comparing estimators $\hat\beta$ according to the criterion of Equation~(\ref{eq:riskdef}) in a \emph{non-asymptotic} fixed $n$, fixed $p>n$ setting. Our approach follows the paradigm of \citet{stein1956}, with risk results that are valid for all $p$ and $n$, rather than asymptotic oracle properties. Our specific contribution, established by Theorem~\ref{th:riskgl} and Corollary~\ref{cor:gl}, is to identify the shortcomings of some commonly used global shrinkage estimators in prediction, with shrinkage driven by a single tuning parameter, and to demonstrate that under certain conditions, suitably-chosen component-specific local shrinkage parameters can result in theoretically lower predictive risk.

\subsection{Connections with existing global shrinkage regression approaches}
We define shrinkage estimators with a single tuning parameter as ``global.'' Examples include ridge regression \citep{hoerl70} and principal components regression or PCR \citep{jolliffe82}, and they remain popular in prediction under the high-dimensional model of Equation (\ref{eq:reg1}). Shrinkage methods enjoy a number of advantages over simultaneous shrinkage and selection-based methods such as the lasso \citep{tibshirani96} and comfortably outperform them in predictive performance in certain situations. Prominent among these is when the predictors are correlated and the resulting lasso estimate is unstable, whereas ridge or PCR estimates are not \citep[see, e.g, the discussion in Chapter 3 of][]{hastie09}. On the theoretical side, \citet{polson12} used a representation devised by  \citet{frank93} to show that many commonly used high-dimensional shrinkage regression estimates, such as the estimates of ridge regression, regression with g-prior \citep{zellner86} and PCR, can be viewed as posterior means under a unified framework of ``global'' shrinkage prior on the regression coefficients that are suitably orthogonalized.  \citet{polson12} also demonstrated that purely global shrinkage regression methods suffer from two major difficulties: (i) the amount of relative shrinkage is monotone in the singular values of the design matrix and (ii) the shrinkage is determined by a single tuning parameter. Both of these factors can translate to poor out of sample prediction performance, which they demonstrated numerically.

\citet{polson12} further provided numerical evidence that the difficulties mentioned above can be resolved by allowing ``local,'' component-specific shrinkage terms, in conjunction with a global shrinkage parameter as used in ridge or PCR, giving rise to the so-called ``global-local'' shrinkage regression models. Specifically, \citet{polson12} demonstrated by simulations that using the horseshoe prior of  \citet{carvalho2010horseshoe} on the regression coefficients performed well over a variety of competitors in terms of predictive performance, including the lasso, ridge, PCR and sparse partial least squares \citep{chun2010sparse}. While these empirical results are encouraging, a theoretical investigation of the conditions required for the horseshoe regression to outperform a global shrinkage regression method such as ridge or PCR in terms of prediction has been lacking. Our work bridges this theoretical gap by developing formal tools for comparing the finite sample predictive risk for shrinkage methods.

\subsection{Regression with non-convex penalties} The $\ell_1$ or $\ell_2$ penalties that correspond to lasso or ridge regression are convex. While this simplifies the computation, it also results in a number of drawbacks such as bias in estimating large signals \citep{fan2001variable}. This problem can be remedied using non-convex $\ell_q$ penalties for $0 <q<1$, but this introduces other problems such as non-uniqueness of solutions and greater computational burden. Prominent examples of non-convex penalties include the smoothly clipped absolute deviation or SCAD  \citep{fan2001variable} and the minimax concave penalty of MCP \citep{zhang2010nearly}. In particular, the MCP estimate enjoys a number of asymptotic optimality properties and conditions required for an iterative computational algorithm to reach the global optimum are available \citep{mazumder2011sparsenet}. The optimality results, however, are valid only in an asymptotic regime with various assumptions on the design $X$ and do not characterize finite sample risk properties. Nevertheless, the univariate MCP estimator is identical to the firm shrinkage estimator of \citet{gao1997waveshrink}, who provide explicit finite sample expressions for predictive risk, a fact utilized later in Section~\ref{sec:comparemcp}.

\subsection{Finite sample estimates of predictive risk} The quadratic risk in Equation~(\ref{eq:riskdef}) involves the future observation $y^*$ and must be estimated. Developing a formal estimate based on the training data $(X,y)$ to compare predictive performance of competing regression methods is important in both frequentist and Bayesian settings. This is because the frequentist tuning parameter or the Bayesian hyper-parameters can then be chosen to minimize the estimated predictive risk, if prediction of future observations is the main modeling goal. A \emph{finite sample unbiased estimate} of $R$ in Equation (\ref{eq:riskdef}) is given by
Stein's unbiased risk estimate or SURE \citep{stein81}.

We will focus on SURE as our estimate of $R$ for the remainder of this article, which is an example of a model-based covariance penalty. Other examples of covariance penalties include Mallows' $C_p$ \citep{mallows73}, Akaike's information criterion \citep{akaike74} and risk inflation criterion \citep{foster1994risk}. Nonparametric penalties include the generalized cross validation of \citet{craven78}, which has the advantage of being model free but usually produces a prediction error estimate with high variance \citep{efron83}. The relationship between the covariance penalties and nonparametric approaches were further explored by \citet{efron04}, who showed the covariance penalties to be a Rao-Blackwellized version of the nonparametric penalties. Thus, \citet{efron04} concluded that model-based penalties such as SURE or Mallows' $C_p$ (the two coincide for models where the fit is linear in the response variable) offer substantially lower variance in estimating the prediction risk, assuming of course the model is true. From a computational perspective, calculating SURE, when it is explicitly available, is substantially less burdensome than performing cross validation, which usually requires several Monte Carlo replications. Furthermore, SURE, an estimate of quadratic risk in prediction, also has connections with the Kullback--Leiber risk for the predictive density \citep{george2006}. 

\subsection{Outline of main contributions} 
Our main contribution is to analyze the finite sample predictive risk of global shrinkage regression methods, examine where these methods fall short, and demonstrate a remedy using local shrinkage parameters. The main  results are summarized as follows.

\subsubsection{Theoretical findings} The key technique to our innovation is an orthogonalized representation first employed by \citet{frank93} that allows shrinkage regression estimates to be viewed as posterior means under some suitable priors. This is formulated in Section~\ref{sec:shrink-postmean}. Using this representation in Sections~\ref{sec:global} and~\ref{sec:g-l}, we devise general, explicit and numerically stable techniques for computing SURE for regression models that can be employed to compare the performances of global as well as horseshoe regressions. We characterize the finite sample risk properties of all competing methods by computing expectations of SURE. Consequently, all results provided in our article are valid under minimal assumptions on the design matrix, similar to the risk results by \citet{stein1956}, where the only requirement is $n>2$. This is at a contrast with most existing results in linear regression focusing on asymptotic minimax risk that require various assumptions on the singular values of $X$ \citep[e.g.,][]{raskutti2011minimax, castillo2015bayesian, dobriban2015high}.

Using the developed tools for SURE, we provide explicit finite sample risk comparisons between the global ridge and global-local horseshoe regressions in Section~\ref{sec:risk}, where for all methods the tuning parameter is chosen to optimize SURE. Specifically, we demonstrate that the horseshoe regression can outperform the optimal ridge regression in prediction when most true signals are zero, but a few are large. We also compare risk of the horseshoe regression with non-convex penalized likelihood approaches such as MCP in Section~\ref{sec:comparemcp} and show that when most of the true signals are away from zero, the risk of MCP can be quite large, unlike that of the horseshoe regression.


\subsubsection{Empirical findings} Extensive numerical results are provided in Section~\ref{sec:sim} and Supplementary Section~\ref{sec:suppsim}. Our simulation results treat three distinct regimes: (i) sparse-robust: where most true signals are zero and a few are large, (ii) null: where all signals are zero and (iii) dense: where all signals are large. Our major finding is that the horseshoe regression outperforms the other methods in (i). Moreover, it is not much worse than ridge in (iii) and adaptive lasso in (ii), which are usually the best performers in these settings. Being a shrinkage estimate, the results for the horseshoe are numerically stable, unlike that of the selection-based estimators in the dense case. We conclude with a demonstration on real data in Section~\ref{sec:real} and by outlining some possible extensions of the current work in Section~\ref{sec:conc}.

\section{Shrinkage regression estimates as posterior means}\label{sec:shrink-postmean}

Let $X=UDW^{T}$ be the singular value decomposition of the design matrix. Let $D=\mathrm{diag}(d_1, \ldots, d_n)$ with $d_1 \geq \ldots \geq d_n >0$ and $\mathrm{Rank}(D) = \min(n,p)=n$. Define $Z=UD$ and  $\alpha=W^{T}\beta$. Then the regression model of Equation (\ref{eq:reg1}) can be reformulated as:
\begin{eqnarray}
y=Z \alpha + \epsilon.\label{eq:reg2}
\end{eqnarray}
The ordinary least squared (OLS) estimate of $\alpha$ is $\hat \alpha = (Z^T Z)^{-1} Z^{T}y= D^{-1} U^{T} y$. Following the original results by \citet{frank93}, several authors have used the well-known orthogonalization technique \citep{polson12, clyde1996prediction, denison12} to demonstrate that the estimates of many shrinkage regression methods can be expressed in terms of the posterior mean of the ``orthogonalized'' regression coefficients $\alpha$ under the following hierarchical model:
\begin{eqnarray}
(\hat \alpha_i \mid \alpha_i, \sigma^2)  &\stackrel{ind}\sim& \Nor(\alpha_i, \sigma^2 d_i^{-2}), \label{eq:alphahat} \\
(\alpha_i \mid \sigma^2, \tau^2,  \lambda_i^2) &\stackrel{ind}\sim& \Nor (0, \sigma^2 \tau^2 \lambda_i^2), \label{eq:tau}
\end{eqnarray}
with $\sigma^2, \tau^2>0$. The global term $\tau$ controls the amount of shrinkage and the fixed $\lambda_i^2$ terms depend on the method at hand. Given $\lambda_i$ and $\tau$, the estimate for $\beta$ under the global shrinkage prior, denoted by $\tilde \beta$, can be expressed in terms of the posterior mean estimate for $\alpha$ as follows:
\begin{eqnarray}
\tilde \alpha_i = \frac{\tau^2 \lambda_i^2 d_i^2}{1 + \tau^2 \lambda_i^2 d_i^2} \hat \alpha_i, \quad \text{and}\quad \tilde \beta = \sum_{i=1}^{n} \tilde \alpha_i w_i,\label{eq:shrink}
\end{eqnarray}
where $\tilde \alpha_i = E (\alpha_i \mid \tau, \lambda_i^2, X, y)$; $w_i$ is a $p \times 1$ vector and is the $i$th column of the $p \times n$ matrix $W$ and the term $\tau^2 \lambda_i^2 d_i^2/(1 + \tau^2 \lambda_i^2 d_i^2) \in (0,1)$ is the shrinkage factor. The expression from Equation (\ref{eq:shrink}) makes it clear that it is the orthogonalized OLS estimates $\hat\alpha_i$s that are shrunk. We shall show that this orthogonalized representation is also particularly suitable for calculating the prediction risk estimate. The reason is tied to the independence assumption that is now feasible in Equations (\ref{eq:alphahat}) and (\ref{eq:tau}). To give a few concrete examples, we note below that several popular shrinkage regression models fall under the framework of Equations (\ref{eq:alphahat}--\ref{eq:tau}):
\begin{enumerate}
\item For ridge regression, $\lambda_i^2=1,  \forall i$, and we have $\tilde \alpha_i=\tau^2 d_i^2\hat \alpha_i/(1+ \tau^2 d_i^2)$.
\item For $K$ component PCR, $\lambda_i^2$ is infinite for the first $K$ components and then 0. Thus, $\tilde \alpha_i= \hat\alpha_i $ for $i=1, \ldots, K$ and $\tilde \alpha_i = 0$ for $i=K+1, \ldots, n$.
\item For regression with g-prior, $\lambda_i^2 = d_i^{-2}$ and we have $\tilde \alpha_i = \tau^2 \hat \alpha_i /(1+\tau^2)$ for $i=1, \ldots, n$.
\end{enumerate}
This shows the amount of relative shrinkage $\tilde \alpha_i/\hat \alpha_i$ is constant in $d_i$ for PCR and g-prior and is monotone in $d_i$ for ridge regression. In none of these cases it depends on the OLS estimate $\hat \alpha_i$ (consequently, on $y$). In the next section we quantify the effect of this behavior on the prediction risk estimate.

\section{Stein's unbiased risk estimate for global shrinkage regression}\label{sec:global}
Define the fit $\tilde y = X\tilde \beta = Z\tilde \alpha$, where $\tilde \alpha$ is the posterior mean of $\alpha$. As noted by \citet{stein81}, the fitted risk is an underestimation of the prediction risk, and SURE for prediction is defined as:
$$
SURE =  ||y- \tilde y||^2 + 2\sigma^2 \sum_{i=1}^{n} \frac{\partial \tilde y_i}{\partial y_i},
$$
where the $\sum_{i=1}^{n} ({\partial \tilde y_i}/{\partial y_i})$ term is also known as the ``degrees of freedom'' \citep{efron04}. By Tweedie's formula \citep{masreliez75, pericchi1992exact} that relates the posterior mean with the marginals; we have for a Gaussian model of Equations (\ref{eq:alphahat}--\ref{eq:tau}) that: $\tilde \alpha = \hat \alpha + \sigma^2 D^{-2} \nabla_{\hat \alpha} \log m(\hat \alpha)$, where $m(\hat\alpha)$ is the marginal for $\hat \alpha$.  Noting $y=Z\hat \alpha$  yields $ \tilde y = y + \sigma^2 UD^{-1} \nabla_{\hat \alpha} \log m(\hat \alpha)$. Using the independence of $\alpha_i$s, the formula for SURE becomes 
\begin{align}
SURE =&  \sigma^4 \sum_{i=1}^{n} d_i^{-2} \left \{ \frac{\partial}{\partial \hat \alpha_i} \log m(\hat \alpha_i)\right \}^2 + 2\sigma^2 \sum_{i=1}^{n} \left\{1 + \sigma^2d_i^{-2} \frac{\partial^2}{\partial \hat \alpha_i^2 } \log m(\hat \alpha_i)\right\}.\label{eq:sure}
\end{align}
Thus, the prediction risk estimate for shrinkage regression can be quantified in terms of the first two derivatives of the log marginal for $\hat\alpha$. Integrating out $\alpha_i$ from Equations (\ref{eq:alphahat}--\ref{eq:tau}) yields in all these cases,
$$
(\hat\alpha_i \mid \sigma^2, \tau^2, \lambda_i^2) \stackrel{ind}\sim \Nor (0, \sigma^2(d_i^{-2} + \tau^2 \lambda_i^2)).
$$
The marginal of $\hat \alpha$ is given by
$$
m(\hat\alpha)  \propto  \prod_{i=1}^{n} \exp \left \{-\frac{\hat \alpha_i^2/2}{\sigma^2(d_i^{-2} + \tau^2 \lambda_i^2)} \right\},
$$
which yields
\begin{eqnarray}
\quad \;  \frac{\partial \log m(\hat\alpha_i)}{\partial \hat\alpha_i} = \frac{-\hat \alpha_i}{\sigma^2(d_i^{-2} + \tau^2 \lambda_i^2)}; \quad \quad \frac{\partial^2 \log m(\hat\alpha_i)}{\partial \hat\alpha_i^2} =  \frac{-1}{\sigma^2(d_i^{-2} + \tau^2 \lambda_i^2)}.\label{eq:logm}
\end{eqnarray}
Therefore, Equation (\ref{eq:sure}) reduces to the following expression for SURE for global shrinkage regressions: $SURE=\sum_{i=1}^{n} SURE_i$, where,
\begin{eqnarray}
SURE_i &=&   \frac{\hat \alpha_i^2 d_i^2}{(1+ \tau^2 \lambda_i^2 d_i^2)^2} + 2\sigma^2 \frac{\tau^2 \lambda_i^2 d_i^2}{(1 + \tau^2 \lambda_i^2 d_i^2)}.\label{eq:risk}
\end{eqnarray}
From a computational perspective, the expression in Equation (\ref{eq:risk}) is attractive, as it avoids costly matrix inversions. For a given $\sigma$ one can choose $\tau$ to minimize the prediction risk, which amounts to a one-dimensional optimization. Note that in our notation, $d_1 \geq d_2 \ldots \geq d_n >0$. Clearly, this is the SURE when $\lambda_i$s are fixed and finite (e.g., ridge regression). For $K$ component PCR, only the first $K$ terms appear in the sum.  The $d_i$ terms are features of the design matrix $X$ and one may try to control the prediction risk by varying $\tau$. When $\tau\to\infty$, $SURE \to 2n\sigma^2$, the risk of prediction with ordinary least squares (unbiased). When $\tau\to0$, we get the mean-only model (zero variance): $ SURE \to  \sum_{i=1}^{n} \hat\alpha_i^2 d_i^2.$
Regression models with $\tau \in (0, \infty)$ represent a bias-variance tradeoff. Following are the two major difficulties of global shrinkage regression:
\begin{enumerate}
\item The first term of Equation~(\ref{eq:risk}) shows that $SURE$ is increased by those components for which $\hat \alpha_i^2 d_i^2$ is large. Choosing a large $\tau$ alleviates this problem, but at the expense of an $SURE_i$ of $2\sigma^2$ even for components for which $\hat \alpha_i^2 d_i^2$ is small (due to the second term in Equation (\ref{eq:risk})). Thus, it might be beneficial to \emph{differentially minimize} the effect of the components for which $\hat \alpha_i^2 d_i^2$ is large, while ensuring those for which $\hat \alpha_i^2 d_i^2$ is small make a contribution less than $2\sigma^2$ to $SURE$. 
Yet, regression models with $\lambda_i$ fixed, such as ridge, PCR, regression with g-priors, provide no mechanism for achieving this, since the relative shrinkage, defined as the ratio $\tilde \alpha_i/\hat\alpha_i$, equals $\tau^2 \lambda_i^2 d_i^2/(1 + \tau^2 \lambda_i^2 d_i^2)$, and is solely driven by a single quantity $\tau$. 
 
\item  Equation~(\ref{eq:shrink}) shows that the relative shrinkage for $\hat\alpha_i$ is monotone in $d_i$; that is, those $\hat\alpha_i$ corresponding to a smaller $d_i$ are necessarily shrunk more (in a relative amount). This is only sensible in the case where one has reasons to believe the low variance eigen-directions (i.e., principal components) of the design matrix are not important predictors of the response variables, an assumption that can be violated in real data \citep{polson12}. 
\end{enumerate}
In the light of these two problems, we proceed to demonstrate that putting a heavy-tailed prior on $\lambda_i$s, in combination with a suitably small value of $\tau$ to enable global-local shrinkage can resolve both these issues. The intuition behind this is that a small value of a \emph{global} parameter $\tau$ enables shrinkage towards zero for all the components while the heavy tails of the \emph{local} or component-specific $\lambda_i$ terms ensure the components with large values of $\hat\alpha_i d_i$ are not shrunk too much, and allow the $\lambda_i$ terms to be learned from the data. Simultaneously ensuring both of these factors helps in controlling the prediction risk for both the noise as well as the signal terms.

\section{Stein's unbiased risk estimate for the horseshoe regression}\label{sec:g-l}
The global-local horseshoe shrinkage regression of \citet{polson12} extends the global shrinkage regression models of the previous section by putting a local (component-specific), heavy-tailed half-Cauchy prior on the $\lambda_i$ terms that allow these terms to be learned from the data, in addition to a global $\tau$. The model equations become:
\begin{eqnarray}
(\hat \alpha_i \mid \alpha_i, \sigma^2)  &\stackrel{ind}\sim& \Nor (\alpha_i, \sigma^2 d_i^{-2}), \label{eq:alphahat1} \\
(\alpha_i \mid \sigma^2, \tau^2,  \lambda_i^2) &\stackrel{ind}\sim& \Nor (0, \sigma^2 \tau^2 \lambda_i^2), \label{eq:alpha1} \\
\lambda_i &\stackrel{ind}\sim& C^{+} (0,1), \label{eq:tau1}
\end{eqnarray}
with $\sigma^2, \tau^2 >0$ and $C^{+} (0,1)$ denotes a standard half-Cauchy random variable with density $p(\lambda_i) = (2/\pi)(1+\lambda_i^2)^{-1}$. The posterior mean $\tilde \alpha$ and the regression estimate $\tilde \beta$ are then obtained analogously to Equation (\ref{eq:shrink}), with the only difference being one uses the posterior mean $\E(\lambda_i\mid \hat\alpha_i, \tau)$ instead of a fixed $\lambda_i$.  The marginal prior on $\alpha_i$s that is obtained as a normal scale mixture by integrating out $\lambda_i$s from Equations (\ref{eq:alpha1}) and (\ref{eq:tau1}) is called the horseshoe prior \citep{carvalho2010horseshoe}. Improved mean square error over competing approaches in regression has been empirically observed by \citet{polson12} with horseshoe prior on $\alpha_i$s. The intuitive explanation for this improved performance is that a heavy tailed prior of $\lambda_i$ leaves the large $\alpha_i$ terms of Equation (\ref{eq:alpha1}) un-shrunk in the posterior, whereas the global $\tau$ term provides shrinkage towards zero for all components \citep[see, for example, the discussion by][and the references therein]{polsonscott2012, bhadra2015horseshoe+,carvalho2010horseshoe}. However, no explicit formulation of the prediction risk under horseshoe shrinkage is available so far and we demonstrate below the heavy-tailed priors on $\lambda_i$ terms, in addition to a global $\tau$, can be beneficial in controlling the overall prediction risk. 

Under the model of Equations (\ref{eq:alphahat1}--\ref{eq:tau1}), after integrating out $\alpha_i$ from the first two equations, we have,
$$
(\hat\alpha_i \mid \sigma^2, \tau^2, \lambda_i^2) \stackrel{ind}\sim \Nor (0, \sigma^2(d_i^{-2} + \tau^2 \lambda_i^2)).
$$
We have, $p(\lambda_i) \propto 1/(1+\lambda_i^2)$. Thus, the marginal of $\hat\alpha$, denoted by $m(\hat\alpha)$, is given up to a constant of proportionality by
\begin{align}
m(\hat \alpha) =& \prod_{i=1}^{n} \int_{0}^{\infty} \Nor (\hat \alpha_i \mid 0, \sigma^2(d_i^{-2} + \tau^2 \lambda_i^2) )p (\lambda_i) d\lambda_i \nonumber\\
\propto& (2\pi\sigma^2)^{-n/2}  \prod_{i=1}^{n} \int_{0}^{\infty} \exp \left \{-\frac{\hat \alpha_i^2 d_i^2 /2}{\sigma^2( 1 +  \tau^2d_i^2 \lambda_i^2)} \right\} \frac{d_i }{( 1 + \tau^2 d_i^2 \lambda_i^2)^{1/2}}  \frac{1}{1 + \lambda_i^2} d\lambda_i.\label{eq:m}
\end{align}
This integral involves the normalizing constant of a compound confluent hypergeometric distribution that can be computed using a result of \citet{gordy98}.
\begin{proposition} \citep{gordy98}. \label{prop:gordy}
The compound confluent hypergeometric (CCH) density is given by
\begin{eqnarray*}
\mathrm{CCH}(x; p, q, r, s, \nu, \theta) &=& \frac{x^{p-1} (1-\nu x)^{q-1} \{\theta+ (1-\theta) \nu x\}^{-r} \exp(-sx)}{B(p,q) H(p,q,r,s,\nu, \theta)},
\end{eqnarray*}
for $0<x<1/\nu$, where the parameters satisfy $p>0$, $q>0$, $r \in \mathbb{R}$, $s \in \mathbb{R}$, $0\leq \nu \leq 1$ and $\theta>0$. Here $B(p,q)$ is the beta function and the function $H(\cdot)$ is given by
\begin{eqnarray*}
H(p,q,r,s,\nu,\theta)=\nu^{-p} \exp(-s/\nu)\Phi_1 (q, r, p+q, s/\nu, 1-\theta),
\end{eqnarray*}
where $\Phi_1$ is the confluent hypergeometric function of two variables, given by
\begin{eqnarray}
\Phi_1 (\alpha, \beta, \gamma, x_1, x_2)&=& \sum_{m=0}^{\infty} \sum_{n=0}^{\infty} \frac{(\alpha)_m (\beta)_n}{(\gamma)_{m+n} m! n!} x_1^m x_2^n,\label{eq:Phi}
\end{eqnarray}
where $(a)_k$ denotes the rising factorial with $(a)_0=1, (a)_1=a$ and $(a)_k = (a+k-1)(a)_{k-1}$.
\end{proposition}
We present our first result in the following theorem and show that the marginal $m(\hat\alpha)$ and all its derivatives lend themselves to a series representation in terms of the first and second moments of a random variable that follows a CCH distribution. Consequently, we quantify SURE for the horseshoe regression.
\begin{theorem}\label{th:sure}
Denote $m'(\hat \alpha_i) = (\partial/\partial \hat \alpha_i) m (\hat \alpha_i)$ and $m''(\hat \alpha_i) = (\partial^2/\partial \hat \alpha_i^2) m (\hat \alpha_i)$. Then, the following holds.
\begin{enumerate}[label= \Alph*.] 
\item SURE for the horseshoe shrinkage regression model defined by Equations (\ref{eq:alphahat1}--\ref{eq:tau1}) is given by $SURE= \sum_{i=1}^{n} {SURE}_i$, where the component-wise contribution $SURE_i$ is given by
\begin{eqnarray}
SURE_i &=& 2\sigma^2  - \sigma^4 d_i^{-2} \left\{\frac{m'(\hat\alpha_i)}{m(\hat\alpha_i)}\right\}^2 + 2\sigma^4 d_i^{-2} \frac{m''(\hat\alpha_i)}{m(\hat\alpha_i)}. \label{eq:stein}
\end{eqnarray}
\item Under independent standard half-Cauchy prior on $\lambda_i$s, for the second and third terms in Equation (\ref{eq:stein}) we have:
\begin{equation*}
\frac{m'(\hat\alpha_i)}{m(\hat\alpha_i)} = -\frac{\hat{\alpha}_i d_i^2}{\sigma^2}\E (Z_i),
\quad \text{and} \quad
\frac{m''(\hat\alpha_i)}{m(\hat\alpha_i)} = -\frac{d_i^2}{\sigma^2}\E (Z_i) +\frac{\hat{\alpha}_i^2 d_i^4}{\sigma^4}\E (Z_i^2),
\end{equation*}
where $(Z_i \mid \hat\alpha_i, \sigma, \tau)$ follows a $\mathrm{CCH}(p=1,q=1/2,r=1,s=\hat{\alpha}_i^2d_i^2/2\sigma^2,v=1,\theta=1/\tau^2d_i^2)$ distribution.
\end{enumerate} 
\end{theorem}
A proof is given in Appendix~\ref{app1}. Theorem~\ref{th:sure} provides a computationally tractable mechanism for calculating SURE for the horseshoe shrinkage regression in terms of the moments of CCH random variables. \citet{gordy98} provides a simple formula for all integer moments of CCH random variables. Specifically, he shows if $X\sim \mathrm{CCH}(x; p, q, r, s, \nu, \theta)$ then
\begin{eqnarray}
\E(X^k) = \frac{(p)_k}{(p+q)_k} \frac{H(p+k, q,r,s,\nu,\theta)} {H(p, q,r,s,\nu, \theta)},\label{eq:gordymoment}
\end{eqnarray}
for integers $k \ge 1$. Moreover, as demonstrated by \citet{gordy98}, these moments can be numerically evaluated quite easily over a range of parameter values and calculations remain very stable. A consequence of this explicit formula for SURE is that the global shrinkage parameter $\tau$ can now be chosen to minimize SURE by performing a one-dimensional numerical optimization. Another consequence is that an application of Theorem 3 of \citet{carvalho2010horseshoe} shows
$$
\lim_{|\hat\alpha_i|\to \infty}  \frac{m'(\hat\alpha_i)}{m(\hat\alpha_i)} = \lim_{|\hat\alpha_i|\to \infty} \frac{\partial \log m(\hat\alpha_i)}{\partial \hat\alpha_i} = 0,
$$
with high probability, where $m(\hat\alpha_i)$ is the marginal under the horseshoe prior.  Recall that the posterior mean $\tilde \alpha_i$ and the OLS estimate $\hat\alpha_i$ are related by Tweedie's formula as
$
\tilde \alpha_i = \hat\alpha_i + \sigma^2 d_i^{-2} {\partial \log m(\hat\alpha_i)}/{\partial\hat\alpha_i}. 
$
Thus, $\tilde \alpha_i \approx \hat\alpha_i$, with high probability, as $|\hat\alpha_i|\to \infty$, for any fixed $d_i$ and $\sigma$ for the horseshoe regression. Since $\hat\alpha_i$ is unbiased for $\alpha_i$, the resultant horseshoe posterior mean is also seen to be unbiased when $|\hat\alpha_i|$ is large. Compare with the resultant $\tilde \alpha_i$ for global shrinkage regression of Equation (\ref{eq:shrink}), which is monotone decreasing in $d_i$, and therefore can be highly biased if a true large $|\alpha_i|$ corresponds to a small $d_i$. 
Perhaps more importantly, we can use the expression from Theorem~\ref{th:sure} to estimate the prediction risk of the horseshoe regression for the signal and the noise terms. First we treat the case when $|\hat\alpha_i|$ is large. We have the following result.
\begin{theorem} \label{th:large}
Define $s_i=\hat\alpha_i^2d_i^2/ 2\sigma^2$, $\theta_i = (\tau^2 d_i^{2})^{-1}$. For any $s_i \ge 1, \theta_i \ge 1$, we have for the horseshoe regression of Equations (\ref{eq:alphahat1}--\ref{eq:tau1}) that
\small
 \begin{eqnarray*}
\left\{1 -\theta_i(\tilde C_1+ \tilde C_2) \frac{(1+s_i)}{s_i^2}  - \theta_i^2(\tilde C_1+ \tilde C_2)^2\frac{(1+s_i)^2}{s_i^3}  \right\} \le \frac{{SURE}_i}{2\sigma^2} \le \left\{1+ 2 \theta_i (1+s_i) \left( \frac{C_1}{s_i^2} + \frac{C_2}{s_i^{3/2}}\right) \right\},
\end{eqnarray*}
\normalsize
where $C_1 = \{1- 5/(2e)\}^{-{1}/{2}}\approx 3.53, C_2 =16/15,  \tilde C_1 = (1- 2/e)^{-1/2} \approx 1.95, \tilde C_2 =4/3,$ 
are constants.
\end{theorem}
A proof is given in Appendix~\ref{app2}. Our result is non-asymptotic, i.e., it is valid for any $s_i \ge1$. However, an easy consequence is that ${SURE}_i \to 2\sigma^2$, almost surely, as $s_i \to \infty$, provided $\tau^2\le d_i^{-2}$. An intuitive explanation of this result is that component-specific shrinkage is feasible in the horseshoe regression model due to the heavy-tailed $\lambda_i$ terms, which prevents the signal terms from getting shrunk too much and consequently, making a large contribution to SURE due to a large bias. With just a global parameter $\tau$, this component-specific shrinkage is not possible. A comparison of $SURE_i$ resulting from Theorem~\ref{th:large} with that from Equation (\ref{eq:risk}) demonstrates using global-local horseshoe shrinkage, we can rectify a major shortcoming of global shrinkage regression, in that the terms with large $s_i$ do not make a large contribution to the prediction risk. Moreover, the main consequence of Theorem~\ref{th:large}, that is ${SURE}_i \to 2\sigma^2$, almost surely, as $s_i \to \infty$, holds for a larger class of ``global-local'' priors, of which the horseshoe is a special case. 

\begin{theorem} \label{th:gl}
Consider the hierarchy of Equations (\ref{eq:alphahat1}--\ref{eq:alpha1}) and suppose the prior on $\lambda_i$ in Equation (\ref{eq:tau1}) satisfies $p(\lambda_i^2) \sim (\lambda_i^2)^{a-1} L(\lambda_i^2)$ as $\lambda_i^2\to\infty$, where $f(x) \sim g(x)$ means $\lim_{x\to\infty} f(x)/g(x)=1$. Assume $a\le 0$ and $L(\cdot)$ is a slowly-varying function, defined as $\lim_{|x|\to \infty} L(tx)/L(x) =1$ for all $t \in (0,\infty)$. Then we have  ${SURE}_i \to 2\sigma^2$, almost surely, as $s_i \to \infty$. 
\end{theorem}
A proof is given in Appendix~\ref{app:thgl}. Densities that satisfy $p(\lambda_i^2) \sim (\lambda_i^2)^{a-1} L(\lambda_i^2)$ as $\lambda_i^2\to\infty$ are sometimes called regularly varying or heavy-tailed. Clearly, the horseshoe prior is a special case, since for the standard half-Cauchy we have $p(\lambda_i) \propto 1/(1+\lambda_i^2)$, which yields by a change of variables $p(\lambda_i^2) = (\lambda_i^2)^{-3/2} \{\lambda_i^2/(1+\lambda_i^2)\}$, which is of the form  $(\lambda_i^{2})^{a-1} L(\lambda_i^2)$ with $a=-1/2$ since $L(\lambda_i^2) = \lambda_i^2/(1+\lambda_i^2)$ is seen to be slowly-varying. Other priors that fall in this framework are the horseshoe+ prior of \citet{bhadra2015horseshoe+}, for which $p(\lambda_i)\propto \log(\lambda_i)/(\lambda_i^2 - 1) = \lambda_i^{-2} L(\lambda_i^2)$ with $L(\lambda_i^2) = \log(\lambda_i) \lambda_i^2/(\lambda_i^2 -1)$. \citet{ghosh2013asymptotic} show that the generalized double Pareto prior \citep{armagan2013generalized} and the three parameter beta prior \citep{armagan2011generalized} also fall in this framework. Thus, Theorem~\ref{th:gl} generalizes the main consequence of Theorem~\ref{th:large} to a broader class of priors in the asymptotic sense as $s_i\to\infty$.  

Next, for the case when $|\hat\alpha_i|$ is  small, we have the following result for estimating the prediction risk of the horseshoe regression.

\begin{theorem} \label{th:small}
Define $s_i=\hat\alpha_i^2d_i^2/ 2\sigma^2$ and $\theta_i = (\tau^2 d_i^{2})^{-1}$. Then the following statements are true for the horseshoe regression.
\begin{enumerate}[label= \Alph*.]
\item $SURE_i$ is an increasing function of $s_i$ in the interval $s_i \in [0,1]$ for any fixed $\tau$.
\item When $s_i=0$, we have that $SURE_i$ is a monotone increasing function of $\tau$, and is bounded in the interval $(0, 2\sigma^2/3]$ when $\tau^2 d_i^2 \in (0,1]$.
\item When $s_i=1$, we have that $SURE_i$ is bounded in the interval $(0, 1.93\sigma^2]$ when $\tau^2 d_i^2 \in (0,1]$.
\end{enumerate}
\end{theorem}
A proof is given in Appendix~\ref{app3}. This theorem establishes that: (i) the terms with smaller $s_i$ in the interval $[0,1]$ contribute less to SURE, with the minimum achieved at $s_i=0$ (these terms can be thought of as the noise terms) and (ii) if $\tau$ is chosen to be sufficiently small, the terms for which $s_i=0$, has an upper bound on SURE at $2\sigma^2/3$. Note that the OLS estimator has risk $2\sigma^2$ for these terms. At $s_i=0$, the PCR risk is either $0$ or $2\sigma^2$, depending on whether the term is or is not included. A commonly used technique for shrinkage regressions is to choose the global $\tau$ to minimize a data-dependent estimate of the risk, such as $C_L$ or SURE \citep{mallows73}. The ridge regression SURE at $s_i=0$ is an increasing function of $\tau$ and thus, it might make sense to choose a small $\tau$ if all $s_i$ terms were small. However, in the presence of some $s_i$ terms that are large, ridge regression cannot choose a very small $\tau$, since the large $s_i$ terms will then be heavily shrunk and contribute too much to SURE. This is not the case with global-local shrinkage regression methods such as the horseshoe, which can still choose a small $\tau$ to mitigate the contribution from the noise terms and rely on the heavy-tailed $\lambda_i$ terms to ensure large signals are not shrunk too much. Consequently, the ridge regression risk estimate is usually larger than the global-local regression risk estimate even for very small $s_i$ terms, when some terms with large $s_i$ are present along with mostly noise terms. At this point, the results concern the risk estimate (i.e., SURE) rather than risk itself, the discussion of which is deferred until Section~\ref{sec:risk}.


To summarize the theoretical findings, Theorem~\ref{th:large} together with Theorem~\ref{th:small} establishes that the horseshoe regression is effective in handling both very large and very small values of $\hat\alpha_i^2 d_i^2$. Specifically, Theorem~\ref{th:small} asserts that a small enough $\tau$ shrinks the noise terms towards zero, minimizing their contribution to SURE. Whereas, according to Theorem~\ref{th:large}, the heavy tails of the Cauchy priors for the $\lambda_i$ terms ensure the large signals are not shrunk too much and ensures a SURE of $2\sigma^2$ for these terms, which is an improvement over purely global methods of shrinkage. 

\section{Prediction risk for the global and horseshoe regressions}\label{sec:risk} In this section we compare the theoretical prediction risks of global and global-local horseshoe shrinkage regressions. While SURE is a data-dependent estimate of the theoretical risk, these two quantities are equal in expectation for all $n$. We use a concentration argument to derive conditions under which the horseshoe regression will outperform global shrinkage regression, e.g., ridge regression, in terms of predictive risk. While the analysis seems difficult for an arbitrary design matrix $X$, we are able to treat the case of ridge regression for orthogonal design, i.e., $X^TX = I$. Clearly, if the SVD of $X$ is written as $X=UDV^{T}$, then we have $D=I$ and for ridge regression $\lambda_i=1$ for all $i$. Thus, for orthogonal design, Equations (\ref{eq:alphahat}) and (\ref{eq:tau}) become
\begin{eqnarray*}
(\hat \alpha_i \mid \alpha_i, \sigma^2)  &\stackrel{ind}\sim& \Nor(\alpha_i, \sigma^2 ), \\
(\alpha_i \mid \sigma^2, \tau^2,  \lambda_i^2) &\stackrel{ind}\sim& \Nor (0, \sigma^2 \tau^2),
\end{eqnarray*}
where $\tau$ is the global shrinkage parameter. Since the fit in this model is linear in $\hat\alpha_i$, SURE is equivalent to Mallows' $C_L$. Equation (14) of \citet{mallows73} shows that if $\tau$ is chosen to minimize $C_L$, then the optimal ridge estimate is given in closed form by
$$
\alpha_i^\star = \left(1 - \frac{n\sigma^2}{\sum_{i=1}^{n} {\hat\alpha_i^2}}\right)\hat\alpha_i.
$$
Alternatively, the solution can be directly obtained from Equation (\ref{eq:risk}) by taking $d_i=\lambda_i=1$ for all $i$ and by setting $\tau^{\star} = \argmin_\tau \sum_{i=1}^{n} SURE_i$. It is perhaps interesting to note that this ``optimal'' ridge estimate, where the tuning parameter is allowed to depend on the data, is no longer linear in $\hat\alpha$. In fact, the optimal solution $\alpha^\star$ can be seen to be closely related to the James--Stein estimate of $\alpha$ and its risk can therefore be quantified using the risk bounds on the James--Stein estimate. As expected due to the global nature of ridge regression, the relative shrinkage $\alpha_i^\star/\hat\alpha_i$ of the optimal solution only depends on $\vert\hat\alpha\vert^2= \sum_{i=1}^{n} \hat \alpha_i^2$ but not on the individual components of $\hat\alpha$. Theorem 1 of \citet{casella1982limit} shows that 
$$
1 - \frac{n-2}{n+ |\alpha|^2} \le \frac{R(\alpha, \alpha^\star)}{R(\alpha, \hat\alpha)} \le 1 - \frac{(n-2)^2}{n} \left(\frac{1}{n-2 + |\alpha|^2}\right).
$$
Consequently,  if $\vert \alpha \vert^2/n \to c$ as $n\to\infty$ then the James--Stein estimate satisfies
$$
\lim_{n\to\infty} \frac{R(\alpha, \alpha^\star)}{R(\alpha, \hat\alpha)} = \frac{c}{c+1}.
$$
Thus, $\alpha^\star$ offers large benefits over the least squares estimate $\hat\alpha$ for small $c$ but it is practically equivalent to the least squares estimate for large $c$. The prediction risk of the least squares estimate for $p>n$ is simply $2n\sigma^2$, or an average component-specific risk of $2\sigma^2$. We first show that when true $\alpha_i = 0$, the component-specific risk bound of the horseshoe shrinkage regression with a fixed $\tau=1$ (i.e., the case of purely local shrinkage) is less than $2\sigma^2$. We have the following result.
\begin{theorem} \label{th:riskgl} (Prediction risk for the purely local horseshoe regression). 
Let $D= I$ and let the global shrinkage parameter in the horseshoe regression be $\tau^2 = 1$. When true $\alpha_i =0$, an upper bound of the component-wise risk of the purely local horseshoe regression is $1.75\sigma^2< 2\sigma^2$.
\end{theorem}
A proof can be found in Appendix~\ref{app:riskgl}. 
The proof uses the fact that the actual risk can be obtained by computing the expectation of SURE. We split the domains of integration into three distinct regions and use the bounds on SURE from 
Theorems~\ref{th:large} and~\ref{th:small}, as appropriate. 

When true $\alpha_i$ is large enough, a consequence of Theorem~\ref{th:large} is that the component-specific risk for global-local shrinkage regression is $2\sigma^2$. This is because SURE in this case is almost surely equal to $2\sigma^2$ and $\hat\alpha_i$ is concentrated around true $\alpha_i$. Therefore, it is established that if only a few components of true $\alpha$ are large and the rest are zero in such a way that $\vert\alpha\vert^2/n$ is large, then the horseshoe regression with fixed $\tau=1$ outperforms ridge regression in terms of predictive risk. The benefit arises from a lower risk for the $\alpha_i=0$ terms. On the other hand, if all components of true $\alpha$ are zero or all are large, the horseshoe regression need not outperform ridge regression.

Although Theorem \ref{th:riskgl} shows the horseshoe regression with a fixed $\tau=1$ outperforms the optimal ridge regression in predictive risk when $\alpha=0$, a useful corollary is that the optimal horseshoe regression still outperforms the optimal ridge regression, where the optimal global tuning parameters for both methods are chosen by minimizing their respective SURE. 
\begin{Cor}\label{cor:gl} (Prediction risk for the optimal horseshoe regression).
Let $SURE_{HS}(\tau=1)$ and $SURE_{HS}({\tau=\tau_{HS}^{*}})$ denote the SURE for the horseshoe regression with fixed $\tau=1$ and $\tau = \tau_{HS}^{*} = \argmin_\tau {(SURE_{HS}(\tau))}$.  Then, for any $\alpha$, $R(\alpha, \hat\alpha^{HS}({\tau=\tau_{HS}^{*}})) \le  R(\alpha, \hat\alpha^{HS}{(\tau=1)})$.
\end{Cor}
\begin{proof}
\small
\begin{eqnarray*}
R(\alpha, \hat\alpha^{HS}({\tau=\tau_{HS}^{*}})) =  \E_{\hat\alpha|\alpha}(SURE_{HS}(\tau =\tau_{HS}^{*})) \le \E_{\hat\alpha|\alpha}(SURE_{HS}(\tau =1)) =  R(\alpha, \hat\alpha^{HS}{(\tau=1)}). 
\end{eqnarray*}
\normalsize
\end{proof}
Clearly, $\tau_{HS}^{*}$ is a function of the data and this complicates exact prediction risk calculations for the optimal horseshoe regression as an expectation of SURE as in Theorem~\ref{th:riskgl}. It is not clear if an explicit minimizer of SURE analogous to Equation (14) of \citet{mallows73} for ridge regression can be obtained for the horseshoe regression.
Nevertheless, Corollary~\ref{cor:gl} shows the risk for the horseshoe regression can only decrease further if one sets $\tau= \tau^{*}_{HS}$, similar to the risk result of \citet{stein1956}. This holds because the expectations of SURE are computed with respect to the distribution of $\hat\alpha$, which is independent of $\tau$ given true $\alpha$. 

\section{Risk comparisons with other non-convex regressions}\label{sec:comparemcp} 
In this section we compare the risk of the proposed horseshoe regression with other approaches that are not shrinkage methods. Specifically, we consider the minimax concave penalty (MCP) of \citet{zhang2010nearly}. Again, for simplicity assume that the design matrix $X$ is orthogonal. As pointed out by \citet{zhang2010nearly}, in this case the solution to the MCP estimator is available in closed form and reduces to the firm shrinkage estimator of \citet{gao1997waveshrink}, which is given by
\begin{eqnarray*}
\delta_{\lambda, \gamma} (\hat\alpha_i) = 
 \begin{cases}
      0, &  \text{ if } |\hat\alpha_i| \le \lambda,\\
       \mathrm{sign}(\hat\alpha_i) \frac{\gamma(|\hat\alpha_i| - \lambda)}{\gamma - 1},  & \text{ if } \lambda \le |\hat\alpha_i| \le \gamma\lambda,\\
       \hat\alpha_i, & \text{ if } |\hat\alpha_i| >\gamma\lambda,
    \end{cases}
    \end{eqnarray*}
for $\lambda> 0$ and $\gamma>1$. For a fixed $\lambda$, soft and hard thresholding estimators are obtained as $\gamma\to \infty$ and $\gamma\to1+$ respectively. An explicit expression for the risk of this estimator is given in Theorem 1 of \citet{gao1997waveshrink}, from which it can be seen easily that $R(\delta_{\lambda, \gamma}) >\lambda^2 \{1/2 - \Phi(-2\lambda)\}$ when $\alpha_i = \lambda$ for any fixed $\lambda>0$ and $\gamma>1$, where $\Phi(\cdot)$ is the standard normal distribution function. Thus, for MCP to work well, a small value for $\lambda$ is essential. However, $\lambda$ is the threshold below which the estimates are shrunk to zero and a large $\lambda$ is favored in a ``dense" situation, where there are many true parameters several standard deviations away from zero.  While this behavior is not necessarily a problem for the MCP, since it is designed with a sparse situation in mind, it is perhaps desirable to avoid a large risk at a given $\lambda$. The horseshoe regression achieves exactly that, since its component-specific risk in a dense case is $2\sigma^2$ by Theorem~\ref{th:large}. In Supplementary  Section S.1 we verify that the MCP performs worse than both global and global-local shrinkage methods in a dense situation under a variety of designs $X$.


\section{Numerical examples}\label{sec:sim}
We simulate data where $n=100$, and consider the cases $p=100, 200, 300, 400, 500$. Let $B$ be a $p\times k$ factor loading matrix, with all entries equal to 1. Let ${F}_i$ be $k\times 1$ matrix of factor values, with all entries drawn independently from $\Nor (0,1)$. The $i$th row of the $n\times p$ design matrix $X$ is generated by a factor model, with number of factors $k=8$, as follows:
	\begin{equation*}
	{X}_i=BF_i + \xi_i, \quad \xi_i \sim \Nor(0,0.1), \quad \text{for} \quad i=1,\ldots, n.
	\end{equation*}
	Thus, the columns of $X$ are correlated. Let $X=UDW^{T}$ denote the singular value decomposition of $X$. The observations $y$ are generated from Equation (\ref{eq:reg2}) with $\sigma^2=1$, where for the true orthogonalized regression coefficients $\alpha_0$, the 6, 30, 57, 67, and 96th components are randomly selected as signals, and the remaining 95 components are noise terms. Coefficients of the signals are generated by a $\Nor (10, 0.5)$ distribution, and coefficients of the noise terms are generated by a $\Nor (0, 0.5)$ distribution. For the case $n=100$ and $p=500$, some of the true orthogonalized regression coefficients $\alpha_0$, their ordinary least squared estimates $\hat{\alpha}$, and the corresponding singular values $d$ of the design matrix, are shown in Table~\ref{tab:1}. 

        \begin{table}[!t]
    	\centering
    	\caption{The true orthogonalized regression coefficients $\alpha_{0i}$, their ordinary least squared estimates $\hat{\alpha}_i$, and singular values $d_i$ of the design matrix, for $n=100$ and $p=500$.}
    	\label{tab:1}
    	\begin{tabular}{ccccc}
	\toprule
    		$i$ & $\alpha_{0i}$ & $\hat{\alpha}_i$ & $d_i$ & $\hat{\alpha}_i d_i$ \\
    		\toprule   	
    		1 &  0.10   &  0.10   &  635.10 &  62.13 \\    	
    		2 &  -0.44  &  -0.32  &  3.16   &  -1.00 \\
     		\ldots & \ldots & \ldots & \ldots & \ldots \\    	
    		5 &  -0.13  &  0.30  &  3.05   &  0.91 \\
    		6 &  10.07  &  10.22  &  3.02   &  30.88 \\
    		\ldots & \ldots & \ldots & \ldots & \ldots \\
    		29 & 0.46   &  0.60   &  2.53   &  1.53  \\
    		30 & 10.47  &  11.07  &  2.51   &  27.76 \\
    		\ldots & \ldots & \ldots & \ldots & \ldots \\
    		56 & 0.35   &  0.57   &  2.07   &  1.18  \\
    		57 & 10.23  &  10.66  &  2.07   &  22.05 \\
    		\ldots & \ldots & \ldots & \ldots & \ldots \\
    		66 & -0.00  &  -0.35   &  1.90   &  -0.66  \\
    		67 & 11.14  &  11.52  &  1.88   &  21.70 \\
    		\ldots & \ldots & \ldots & \ldots & \ldots \\
    		95 & -0.82  &  -0.56  &  1.42   &  -0.79 \\
    		96 & 9.60   &  10.21   &  1.40   &  14.26 \\
    		\ldots & \ldots & \ldots & \ldots & \ldots \\
    		100 & 0.61  &  0.91   &  1.27   &  1.15  \\
    		\bottomrule
    	\end{tabular}
   \end{table}

 \begin{table}[!t]
    	\centering
    	\caption{SURE and average out of sample prediction SSE (standard deviation of SSE) on one training set and 200 testing sets for the competing methods for $n=100$, for ridge regression (RR), the lasso regression (LASSO), the adaptive lasso (A\_LASSO), principal components regression (PCR) and the horseshoe regression (HS).  The lowest SURE in each row is in italics and the lowest average prediction SSE is in bold. A formula for SURE is unavailable for the adaptive lasso.}
    	\label{tab:result}
	\begin{tabular}{cccccccccc}
    		\toprule
		&\multicolumn{2}{c}{RR} &\multicolumn{2}{c}{LASSO} & A\_LASSO &\multicolumn{2}{c}{PCR} &\multicolumn{2}{c}{HS}\\
    		$p$ & SURE & SSE & SURE & SSE & SSE & SURE & SSE & SURE & SSE \\
    		\toprule	
    		100 & 159.02 & 168.24 & 125.37 & 128.98 & 127.22 & 162.23 & 179.81 & \textit{120.59} & \textbf{126.33} \\    	
    		{} & {} & (23.87) & {} & (18.80) & (18.10) & {} & (25.51) & {} & (18.77) \\
    		200 & 187.38 & 174.92 & 140.99 & 132.46 & 151.89 & 213.90 & 191.33 & \textit{139.32} & \textbf{126.99} \\
    		{} & {} & (21.13) & {} & (18.38) & (20.47) & {} & (22.62) & {} & (17.29) \\ 	
    		300 & 192.78 & 191.91 & \textit{147.83} & 145.04 & 153.64 & 260.65 & 253.00 & 151.24 & \textbf{136.67} \\
    		{} & {} & (22.95) & {} & (19.89) & (21.19) & {} & (26.58) & {} & (18.73) \\
    		400 & 195.02 & 182.55 & 148.56 & 165.63 & 178.98 & 346.19 & 292.02 & \textit{147.69} & \textbf{143.91} \\
    		{} & {} & (22.70) & {} & (21.55) & (20.12) & {} & (28.98) & {} & (18.41) \\
    		500 & 196.11 & 188.78 & 159.95 & \textbf{159.56} & 186.23 & 386.50 & 366.88 & \textit{144.97} & 160.11 \\
    		{} & {} & (22.33) & {} & (19.94) & (23.50) & {} & (39.38) & {} & (20.29) \\
    		\bottomrule
    	\end{tabular}
    \end{table} 
Table~\ref{tab:result} lists the SURE for prediction and actual out of sample sum of squared prediction error (SSE) for the ridge, lasso, PCR and horseshoe regressions. Out of sample prediction error of the adaptive lasso is also included in the comparisons, although we are unaware of a formula for computing the SURE for the adaptive lasso. SURE for ridge and PCR can be computed by an application of Equation (\ref{eq:risk}) and SURE for the horseshoe regression is given by Theorem~\ref{th:sure}. SURE for the lasso is calculated using the result given by \cite{tibshirani2012degrees}.  In each case, the model is trained on 100 samples. We report the SSE on 100 testing samples, averaged over 200 testing data sets, and their standard deviations. For ridge, lasso, PCR and horseshoe regression, the global shrinkage parameters were chosen to minimize SURE for prediction.  In adaptive lasso, the shrinkage parameters were chosen by cross validation due to SURE being unavailable. It can be seen that SURE in most cases are within one standard deviation of the actual out of sample prediction SSE, suggesting SURE is an accurate method for evaluating actual out of sample prediction performance. When $p=100, 200, 300, 400$, horseshoe regression has the lowest prediction SSE. When $p=500$, SSE of the lasso and horseshoe regression are close, and the lasso performs marginally better. The horseshoe regression also has the lowest SURE in all but one cases. Generally, SURE increases with $p$ for all methods. The SURE for ridge regression approaches  the OLS risk, which is $2n\sigma^2=200$ in these situations. SURE for PCR is larger than the OLS risk and PCR happens to be the poorest performer in most settings. Performance of the adaptive lasso also degrades compared to the lasso and the horseshoe, which remain the two best performers. Finally, the horseshoe regression outperforms the lasso in four out of the five settings we considered.

   \begin{figure}[!t]
   	\centering
   	\includegraphics[width=12.5cm, height=8cm]{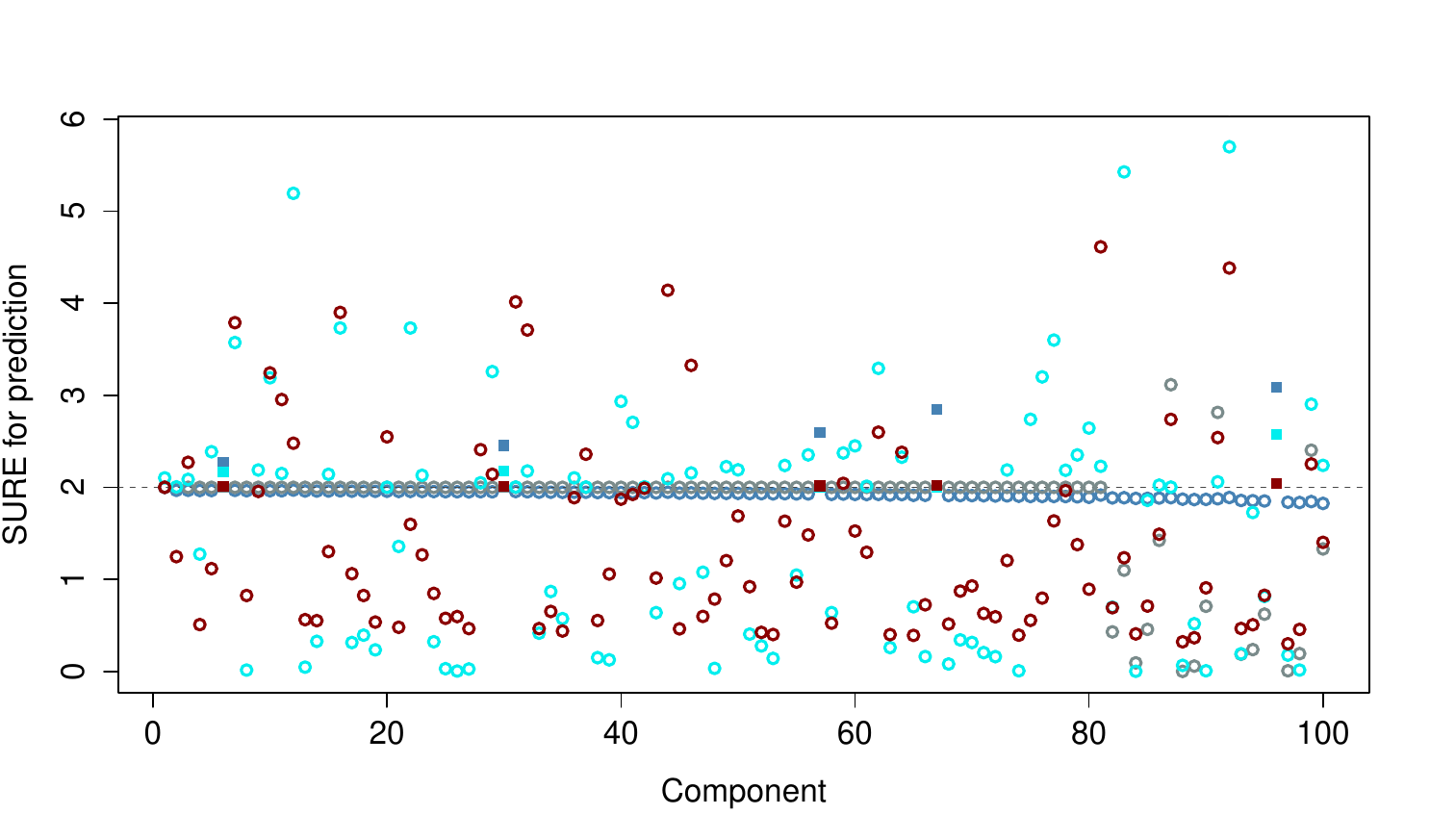}
   	\caption{Component-wise SURE for ridge (blue), PCR (gray), lasso (cyan), and horseshoe regression (red), for $n=100$ and $p=500$. Signal components are shown in solid squares and noise components shown in blank circles. Dashed horizontal line is at $2\sigma^2=2$.}
   	\label{fig:stein}
   \end{figure}
\begin{figure}[!h]
   	\centering
   	\includegraphics[width=12.5cm, height=8cm]{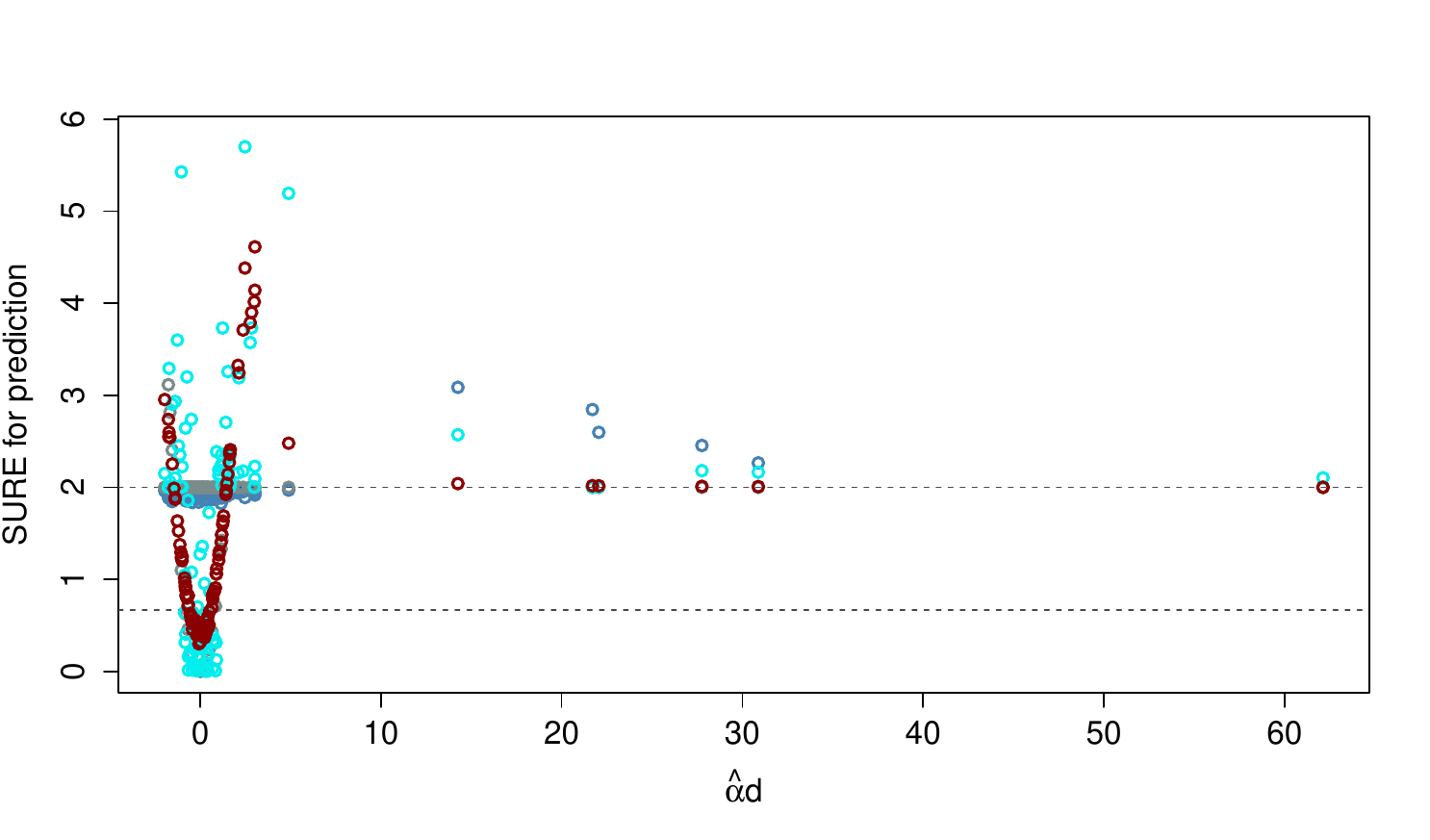}
   	\caption{SURE for ridge (blue), PCR (gray), lasso (cyan) and horseshoe regression (red), versus $\hat{\alpha}d$, where $\hat{\alpha}$ is the OLS estimate of the orthogonalized regression coefficient, and $d$ is the singular value, for $n=100$ and $p=500$. Dashed horizontal lines are at $2\sigma^2=2$ and $2\sigma^2/3 = 0.67$.}
   	\label{fig:2}
   \end{figure}

Figure~\ref{fig:stein} shows contribution to SURE by each component for $n=100$ and $p=500$, for ridge, PCR, lasso and horseshoe regressions. The components are ordered left to right on the $x$-axis by decreasing magnitude of $d_i$, and SURE for prediction on each component are shown on the $y$-axis. Note from Table~\ref{tab:1} that the $6, 30, 57, 67$ and $96$th components are the signals, meaning these terms correspond to a large $\alpha_0$. The PCR risk on the $96$th component is $203.22$, which is out of range for the $y$-axis in the plot.  For this data set, PCR selects 81 components, and therefore SURE for the first $81$ components equal to $2\sigma^2=2$ and the SURE is equal to $\hat\alpha^2_i d_i^2$ for $i=82,\ldots, 100$. Component-wise SURE for ridge regression are large on the signal components, and is decreasing as the singular values $d$ decrease on the other components. But due to the large global shrinkage parameter $\tau$ ridge must select in presence of both large signals and noise terms, the magnitude of improvement over the OLS risk $2\sigma^2$ is small for the noise terms. On the other hand, the horseshoe estimator does not shrink the components with large $\hat\alpha_i d_i$ heavily and therefore the horseshoe SURE on the signal components are almost equal to $2\sigma^2$ (according to Theorem~\ref{th:large}). SURE for the horseshoe is also much smaller than $2\sigma^2$ on many of the noise components. Lasso also appears to be quite effective for the noise terms, but its performance for the signal components is generally not as effective as the horseshoe.

Figure~\ref{fig:2}  takes a fresh look at the same results and shows component-wise SURE plotted against $\hat{\alpha_i}d_i$. The signal components as well as the first component in Table~\ref{tab:1} have $\hat{\alpha_i}d_i>10$. Horseshoe SURE converges to $2\sigma^2$ for large $\hat\alpha_i d_i$, as expected from Theorem~\ref{th:large}. For these components, the SURE for both ridge and lasso are larger than $2\sigma^2$, due to the bias introduced in estimating large signals by these methods \citep[see also Theorem 1 of][]{carvalho2010horseshoe}. When $\hat{\alpha_i}^2d_i^2\approx 0$, risks for lasso and horseshoe are comparable, with lasso being slightly better. This is because an estimate can be exactly zero for the lasso, but not for the horseshoe, which is a shrinkage method (as opposed to a selection method). Nevertheless, the upper bound on SURE for the horseshoe regression at $2\sigma^2/3$ when $\hat{\alpha_i}^2d_i^2\approx 0$ and provided $\tau$ is chosen to be small enough so that $\tau^2 \le d_i^{-2}$, as established by Theorem~\ref{th:small}, can be verified from Figure~\ref{fig:2}.

Additional simulation results are presented in Supplementary Section~\ref{sec:suppsim}, where we (i) treat a higher dimensional case ($p=1000$), (i) perform comparisons with non-convex MCP \citep{zhang2010nearly} and SCAD \citep{fan2001variable} regressions, (iii) explore different choices of $X$ and (iv) explore the effect of the choice of $\alpha$. The main finding is that the horseshoe regression is often the best performer when $\alpha$ has a sparse-robust structure as in Table~\ref{tab:1}, that is most elements are very small while a few are large so that $|\alpha|^2$ is large. This is consistent with the theoretical results of Sections~\ref{sec:risk} and~\ref{sec:comparemcp}. 


\section{Assessing out of sample prediction in a pharmacogenomics data set}\label{sec:real}
We compare the out of sample prediction error of the horseshoe regression with ridge regression, PCR, the lasso, the adaptive  lasso, MCP and SCAD on a pharmacogenomics data set. The data were originally described by \citet{szackas04}, in which the authors studied 60 cancer cell lines in the publicly available NCI-60 database (\href{https://dtp.cancer.gov/discovery\_development/nci-60/}{https://dtp.cancer.gov/discovery\_development/nci-60/}). The goal here is to predict the expression of the human ABC transporter genes (responses) using some compounds or drugs (predictors) at which 50\% inhibition of cellular growth for the cell lines are induced. The NCI-60 database includes the concentration level of 1429 such compounds, out of which we use 853, which did not have any missing values, as predictors. We investigate the expression levels of transporter genes A1 to A12 (except for A11, which we omit due to  missing values), and B1. Thus, in our study $X$ is a $n \times p$ matrix of predictors with $n=60, p=853$ and $Y$ is a $n$-dimensional response vector for each of the 12 candidate transporter genes under consideration.
	
To test the performance of the methods, we split each data set into training and testing sets, with 75\% (45 out of 60) of the observations in the training sets. We standardize each response by subtracting the mean and dividing by the standard deviation. We fit the model on the training data, and then calculate mean squared prediction error (prediction MSE) on the testing data. This is repeated for $20$ random splits of the data into training and testing sets. The tuning parameters in ridge regression, the lasso, the adaptive lasso, SCAD and MCP are chosen by five-fold cross validation on the training data. Similarly, the number of components in PCR and the global shrinkage parameter $\tau$ for horseshoe regression are chosen by cross validation as well. It is possible to use SURE to select the tuning parameters or the number of components, but one needs an estimate of the standard deviation of the errors in high-dimensional regressions. This is a problem of recent interest, as the OLS estimate of $\sigma^2$ is not well-defined in the $p>n$ case. Unfortunately, some of the existing methods we tried, such as the method of moments estimator of \citet{dicker2014variance}, often resulted in unreasonable estimates for $\sigma^2$, such as negative numbers. Thus, we stick to cross validation here, as it is not necessary to estimate the residual standard deviation in that case.   
	
The average prediction MSE over 20 random training-testing splits for the competing methods is reported in Table~\ref{table:3}. Average prediction MSE for  responses A1, A8 and A10 are around or larger than $1$ for all of the methods. Since the responses are standardized before analysis, we might conclude that none of the methods performed well for these cases. 
Among the remaining nine cases, the horseshoe regression substantially outperforms the other methods for A3, A4, A9, A12 and B1. It is comparable to PCR for A5 and A7, and is comparable to the adaptive lasso for A6, which are the best performers in the respective cases. Overall, the horseshoe regression performed the best in 5 among the total 12 cases we considered.
	
	  \begin{table}[!t]
	    	\centering
	    	\caption{Average out of sample mean squared prediction error computed on 20 random training-testing splits (number of splits out of 20 with lowest prediction MSE), for each of the 12 human ABC transporter genes (A1--A10, A12, B1) in the pharmacogenomics example.  Methods under consideration are ridge regression (RR), principal components regression (PCR) , the lasso, the adaptive lasso (A\_LASSO), the minimax concave penalty (MCP), the smoothly clipped absolute deviation (SCAD) penalty, and the horseshoe regression (HS).  Lowest prediction MSE and largest number of splits with the lowest prediction MSE for each response in bold.}
	    	
	       \label{table:3}
	       \begin{tabular}{cccccccc}
	       	\toprule
	       	Response & RR & PCR & LASSO & A\_LASSO & MCP & SCAD & HS \\
	        \toprule    	
	       	A1 &  1.12 & 1.10 & \textbf{1.00} & \textbf{1.00} & 1.01 & 1.06 & 1.30  \\
	       	 & (2) & (5) & \textbf{(7)} & (2) & (1) & (1) & (2) \\    	
	       	A2 &  1.00 & 1.04 & 0.95 & 0.93 & \textbf{0.92} & 0.99 & 1.15  \\
	       	 & (3) & (1) & \textbf{(7)} & (5) & (1) & (0) & (3) \\
	       	A3 &  0.77 & 0.91 & 1.11 & 0.90 & 0.92 & 1.06 & \textbf{0.65}  \\
	       	 & (1) & (0) & (0) & (0) & (1) & (0) & \textbf{(18)} \\
	       	A4 &  0.92 & 0.95 & 0.97 & 0.96 & 0.93 & 0.99 & \textbf{0.79}  \\ 	
	       	 & (2) & (0) & (2) & (2) & (2) & (0) & \textbf{(12)} \\
	       	A5 &  0.82 & \textbf{0.77} & 1.06 & 0.81 & 0.83 & 0.94 & 0.79  \\
	       	 & (1) & \textbf{(6)} & (4) & (1) & (2) & (0) & \textbf{(6)} \\
	       	A6 &  0.93 & 0.92 & 0.98 & \textbf{0.86} & 0.87 & 0.90 & 0.95  \\
	       	 & (4) & (0) & (3) & (5) & (0) & (2) & \textbf{(6)} \\
	       	A7 &  0.92 & \textbf{0.83} & 0.92 & 0.93 & 0.99 & 0.93 & 0.85  \\
	       	 & (0) & \textbf{(8)} & (1) & (4) & (0) & (0) & (7) \\
	       	A8 &  1.08 & 1.05 & 1.14 & \textbf{1.01} & \textbf{1.01} & 1.15 & 1.34  \\
	       	 & \textbf{(6)} & (4) & \textbf{(6)} & (4) & (0) & (0) & (0) \\
	       	A9 &  0.57 & 0.64 & 0.81 & 0.67 & 0.77 & 0.68 & \textbf{0.55}  \\
	       	 & (4) & (0) & (0) & (6) & (0) & (1) & \textbf{(9)} \\
	       	A10 &  1.18 & 1.04 & \textbf{1.00} & 1.01 & \textbf{1.00} & 1.06 & 1.33  \\
	       	 & (0) & \textbf{(7)} & (4) & (3) & (2) & (0) & (4) \\
	       	A12 &  1.01 & 1.12 & 1.09 & 1.01 & 1.02 & 1.05 & \textbf{0.80}  \\
	       	 & (0) & (0) & (2) & (2) & (1) & (0) & \textbf{(15)} \\
	       	B1  &  0.53 & 0.59 & 0.70 & 0.63 & 0.91 & 0.70 & \textbf{0.46}  \\
	       	 & (1) & (0) & (3) & (2) & (1) & (3) & \textbf{(10)} \\
	       	\bottomrule
	       \end{tabular}
	    \end{table}

\section{Concluding remarks}\label{sec:conc}
We outlined some situations where the horseshoe regression is expected to perform better compared to some other commonly used ``global'' shrinkage or selection alternatives for high-dimensional regression. Specifically, we demonstrated that the global term helps in mitigating the prediction risk arising from the noise terms, and an appropriate choice for the tails of the local terms is crucial for controlling the risk due to the signal terms. For this article we have used the horseshoe prior as our choice for the global-local prior. However, in recent years, several other priors have been developed that fall in this class. This includes the horseshoe+ \citep{bhadra2015horseshoe+, bhadra2015default},  the three-parameter beta  \citep{armagan2011generalized}, the normal-exponential-gamma \citep{griffin2005alternative}, the generalized double Pareto \citep{armagan2013generalized}, the generalized shrinkage prior \citep{denison12} and the Dirichlet--Laplace prior \citep{bhattacharya2014dirichlet}. Empirical Bayes approaches have also appeared \citep{martin2014} and the spike and slab priors have made a resurgence due to recently developed efficient computational approaches \citep{rockova2014, rockova2015}. Especially in the light of Theorem~\ref{th:gl}, we expect the results developed in this article for the horseshoe to foreshadow similar results when many of these alternatives are deployed. A particular advantage of using the horseshoe prior seems to be the tractable expression for SURE, as developed in Theorem~\ref{th:sure}. Whether this advantage translates to some of the other global-local priors mentioned above is an open question. 
Following the approach of \citet{stein81}, our risk results are developed in a non-asymptotic setting (finite $n$, finite $p>n$). In the normal means model, finite sample risk properties in estimation under heavy-tailed priors have been considered by \citet{polsonscott2012}. However, their work does not consider (a) predictive risk or (b) a linear regression model. Global-local priors such as the horseshoe and horseshoe+ are known to be minimax in estimation in the Gaussian sequence model \citep{van2014horseshoe, van2016conditions}.  For linear regression, frequentist minimax risk results are discussed by \citet{raskutti2011minimax}; and \citet{castillo2015bayesian} have shown that spike and slab priors achieve minimax prediction risk in regression. Whether the prediction risk for the horseshoe regression is optimal in an asymptotic sense is an important question to investigate and recent asymptotic prediction risk results for ridge regression \citep{dobriban2015high} should prove helpful for comparing with global shrinkage methods. 
Another possible direction for future investigation might be to explore the implications of our findings on the predictive density in terms of an appropriate metric, say the Kullback-Leibler loss, following the results of \citet{george2006}. 

\section*{Acknowledgements}
The authors are grateful for constructive suggestions by the reviewers and the action editor. Bhadra and Polson are supported by Grant No. DMS-1613063 by the US National Science Foundation.

\appendix
\section{Proofs}
\setcounter{equation}{0}
\setcounter{table}{0}
\setcounter{section}{0}
\setcounter{figure}{0}
\setcounter{result}{0}
\renewcommand{\theequation}{A.\arabic{equation}}
\renewcommand{\theresult}{A.\arabic{result}}
\renewcommand{\thesubsection}{A.\arabic{subsection}}
\renewcommand{\thelemma}{A.\arabic{lemma}}
\subsection{Proof of Theorem \ref{th:sure}}\label{app1}
Part A follows from Equation (\ref{eq:sure}) with standard algebraic manipulations. To prove part B, define $Z_i = 1/(1 + \tau^2 \lambda_i^2 d_i^2)$. Then, from Equation (\ref{eq:m})
\begin{eqnarray*}
m(\hat \alpha) & =& (2\pi\sigma^2)^{-n/2}  \prod_{i=1}^{n} \int_{0}^{1}  \exp(-z_i\hat \alpha_i^2 d_i^2/2\sigma^2) d_i z_i^{1/2} \left(\frac{z_i\tau^2d_i^2}{1-z_i+z_i\tau^2d_i^2} \right) \frac{1}{\tau d_i} (1-z_i)^{-1/2}z_i^{-3/2} d z_i \\
&=&  (2\pi\sigma^2)^{-n/2} \prod_{i=1}^{n} \int_{0}^{1}  \exp(-z_i\hat \alpha_i^2 d_i^2/2\sigma^2) (1 - z_i)^{-1/2}\left\{ \frac{1}{\tau^2d_i^2} + \left(1 - \frac{1}{\tau^2 d_i^2}\right) z_i\right\}^{-1} d z_i.
\end{eqnarray*}
From the definition of the compound confluent hypergeometric (CCH) density in \citet{gordy98}, the result of the integral is proportional to the normalizing constant of the CCH density and we have from Proposition~\ref{prop:gordy} that,
\begin{eqnarray*}
m(\hat\alpha) &\propto& (2\pi\sigma^2)^{-n/2} \prod_{i=1}^{n}  H\left( 1, \frac{1}{2} , 1, \frac{\hat \alpha_i^2 d_i^2}{2\sigma^2} , 1, \frac{1}{\tau^2 d_i^2}\right).
\end{eqnarray*}
In addition, the random variable $(Z_i \mid \hat\alpha_i, \sigma, \tau)$ follows a $\mathrm{CCH}(1,1/2,1,\hat{\alpha}_i^2 d_i^2/2\sigma^2,1,1/\tau^2d_i^2)$ distribution. 
Lemma 3 of \citet{gordy98} gives,
$$
\frac{d}{ds} H(p,q,r,s,\nu,\theta) = -\frac{p}{p+q} H(p+1, q, r, s, \nu, \theta).
$$
This yields after some algebra that,
\begin{align*}
\frac{m'(\hat\alpha_i)}{m(\hat\alpha_i)} =& -\frac{2}{3} \frac{H\left(2, \frac{1}{2} , 1, \frac{\hat \alpha_i^2 d_i^2}{2\sigma^2} , 1, \frac{1}{\tau^2 d_i^2}\right)}{H\left( 1, \frac{1}{2} , 1, \frac{\hat \alpha_i^2 d_i^2}{2\sigma^2} , 1, \frac{1}{\tau^2 d_i^2}\right)} \frac{\hat\alpha_i d_i^2}{\sigma^2},\\
\frac{m''(\hat\alpha_i)}{m(\hat\alpha_i)} =& \frac{-\frac{2}{3} H\left(2, \frac{1}{2} , 1, \frac{\hat \alpha_i^2 d_i^2}{2\sigma^2} , 1, \frac{1}{\tau^2 d_i^2}\right)\frac{d_i^2}{\sigma^2} + \frac{8}{15} H\left(3, \frac{1}{2} , 1, \frac{\hat \alpha_i^2 d_i^2}{2\sigma^2} , 1, \frac{1}{\tau^2 d_i^2}\right)\frac{\hat\alpha_i^2 d_i^4}{\sigma^4}}{ H\left( 1, \frac{1}{2} , 1, \frac{\hat \alpha_i^2 d_i^2}{2\sigma^2} , 1, \frac{1}{\tau^2 d_i^2}\right)} .
\end{align*}
The correctness of the assertion
\begin{equation*}
\frac{m'(\hat\alpha_i)}{m(\hat\alpha_i)} = -\frac{\hat{\alpha}_i d_i^2}{\sigma^2}\E (Z_i),
\quad \text{and} \quad
\frac{m''(\hat\alpha_i)}{m(\hat\alpha_i)} = -\frac{d_i^2}{\sigma^2}\E (Z_i)+\frac{\hat{\alpha}_i^2 d_i^4}{\sigma^4}\E (Z_i^2),
\end{equation*}
can then be verified using Equation (\ref{eq:gordymoment}), completing the proof. 

\subsection{Proof of Theorem \ref{th:large}}\label{app2}

Define $s_i=\hat{\alpha}_i^2 d_i^2/2\sigma^2$ and $\theta_i= (\tau^2d_i^2)^{-1}$, with$\theta_i\ge1, s_i\ge 1$.  From Theorem~\ref{th:sure}, the component-wise SURE is
\begin{align}
SURE_i =& 2\sigma^2 - 2\sigma^2\E (Z_i)-\hat{\alpha}_i^2 d_i^2\{\E (Z_i)\}^2
+2\hat{\alpha_i}^2 d_i^2\E (Z_i^2) \nonumber \\
	=& 2\sigma^2 [1-\E (Z_i)+2s_i\E (Z_i^2) -s_i\{\E (Z_i)\}^2], \label{eq:apprisk}
\end{align}
Thus,
\begin{eqnarray*}
2\sigma^2 [1-\E (Z_i) - s_i\{\E(Z_i)\}^2] \le {SURE}_i \le 2\sigma^2 [1+ 2s_i\E (Z_i^2)].
\end{eqnarray*}
To find bounds on SURE, we need upper bounds on $\E(Z_i^2)$ and $\E(Z_i)$.
Clearly, $\theta_i^{-1}\le\{\theta_i+(1-\theta_i)z_i\}^{-1}\le 1$, when $\theta_i\ge 1$. Let $a_i=\log(s_i^{5/2})/s_i$. Then $a_i\in [0,5/(2e))$ when $s_i\ge 1$. Now,
\begin{align*}
\E (Z_i^2) =& \frac{\int_0^1 z_i^2(1-z_i)^{-\frac{1}{2}}\{\theta_i+(1-\theta_i)z_i\}^{-1}\exp(-s_i z_i) d z_i}
{\int_0^1 (1-z_i)^{-\frac{1}{2}}\{\theta_i+(1-\theta_i)z_i\}^{-1}\exp(-s_i z_i) d z_i},
\end{align*}
An upper bound to the numerator of $\E(Z_i^2)$ can be found as follows. 
\begin{align*}
&\int_0^1 z_i^2(1-z_i)^{-\frac{1}{2}}\{\theta_i+(1-\theta_i)z_i\}^{-1}\exp(-s_i z_i) d z_i\\
&\le \int_0^1 z_i^2(1-z_i)^{-\frac{1}{2}}\exp(-s_i z_i) d z_i\\
&= \int_0^{a_i} z_i^2(1-z_i)^{-\frac{1}{2}}\exp(-s_i z_i) d z_i + \int_{a_i}^1 z_i^2(1-z_i)^{-\frac{1}{2}}\exp(-s_i z_i) d z_i\\
&\le (1- a_i)^{-\frac{1}{2}} \int_0^{a_i} z_i^2\exp(-s_i z_i) d z_i + \exp(-a_i s_i) \int_{a_i}^1 z_i^2(1-z_i)^{-\frac{1}{2}} d z_i\\
&= (1- a_i)^{-\frac{1}{2}} \frac{2}{s_i^3}\left\{1 - \left(1+a_i s_i + \frac{a_i^2 s_i^2}{2}\right)\exp(-a_i s_i)\right\} + \exp(-a_i s_i) \int_{a_i}^1 z_i^2(1-z_i)^{-\frac{1}{2}} d z_i\\
&\le \{1- 5/(2e)\}^{-\frac{1}{2}} \frac{2}{s_i^3} + \frac{1}{s_i^{5/2}} \int_0^1 z_i^2(1-z_i)^{-\frac{1}{2}} d z_i\\
&=\frac{C_1}{s_i^3} + \frac{C_2}{s_i^{5/2}},
\end{align*}
where $C_1 = \{1- 5/(2e)\}^{-\frac{1}{2}}\approx 3.53$ and $C_2 = \int_0^1 z_i^2(1-z_i)^{-\frac{1}{2}} d z_i = \Gamma(1/2) \Gamma(3)/\Gamma(3.5) =16/15$. Similarly, a lower bound on the denominator of $\E(Z_i^2)$ is
\begin{align*}
&{\int_0^1 (1-z_i)^{-\frac{1}{2}}\{\theta_i+(1-\theta_i)z_i\}^{-1}\exp(-s_i z_i) d z_i}\\
&\ge \theta_i^{-1}\int_{0}^{1} {\exp(-s_i z_i) d z_i}\\
&=\theta_i^{-1}\left\{\frac{1-\exp(-s_i)}{s_i}\right\} \ge \frac{1}{\theta_i (1+s_i)},
\end{align*}
Thus, combining the upper bound on the numerator and the lower bound on the denominator
\begin{align*}
\E(Z_i^2) &\le \theta_i (1+s_i) \left( \frac{C_1}{s_i^3} + \frac{C_2}{s_i^{5/2}}\right).
\end{align*}
Thus,
\begin{eqnarray}
{SURE}_i &\le& 2\sigma^2 [1+ 2s_i\E (Z_i^2)]\nonumber\\
&\le &2\sigma^2 \left\{1+ 2 \theta_i (1+s_i) \left( \frac{C_1}{s_i^2} + \frac{C_2}{s_i^{3/2}}\right) \right\}. \label{eq:aupper}
\end{eqnarray}
An upper bound to the numerator of $\E(Z_i)$ can be found as follows. Let $\tilde a_i = \log(s_i^2)/s_i$. Then, $\tilde a_i \in [0,2/e)$ for $s_i\ge 1$.
\begin{align*}
&\int_0^1 z_i(1-z_i)^{-\frac{1}{2}}\{\theta_i+(1-\theta_i)z_i\}^{-1}\exp(-s_i z_i) d z_i\\
&\le \int_0^1 z_i(1-z_i)^{-\frac{1}{2}}\exp(-s_i z_i) d z_i\\
&= \int_0^{ \tilde a_i} z_i(1-z_i)^{-\frac{1}{2}}\exp(-s_i z_i) d z_i + \int_{ \tilde a_i}^1 z_i(1-z_i)^{-\frac{1}{2}}\exp(-s_i z_i) d z_i\\
&\le (1- \tilde a_i)^{-\frac{1}{2}} \int_0^{\tilde a_i} z_i\exp(-s_i z_i) d z_i + \exp(-\tilde a_i s_i) \int_{a_i}^1 z_i(1-z_i)^{-\frac{1}{2}} d z_i\\
&= (1- \tilde a_i)^{-\frac{1}{2}} \frac{1}{s_i^2}\left\{1 - \left(1+ \tilde a_i s_i \right)\exp(-\tilde a_i s_i)\right\} + \exp(-\tilde a_i s_i) \int_{\tilde a_i}^1 z_i(1-z_i)^{-\frac{1}{2}} d z_i\\
&\le (1- 2/e)^{-\frac{1}{2}} \frac{1}{s_i^2} + \frac{1}{s_i^{2}} \int_0^1 z_i(1-z_i)^{-\frac{1}{2}} d z_i\\
&=\frac{\tilde C_1}{s_i^2} + \frac{\tilde C_2}{s_i^{2}},
\end{align*}
where $\tilde C_1 = (1- 2/e)^{-1/2} \approx 1.95$ and $\tilde C_2 = \int_0^1 z_i(1-z_i)^{-\frac{1}{2}} d z_i\ = \Gamma(1/2) \Gamma(2)/\Gamma(2.5) =4/3$. The lower bound on the denominator is the same as before. 
Thus,
\begin{align*}
\E(Z_i) &\le \frac{\theta_i(1+s_i)}{s_i^2} \left( \tilde C_1+ \tilde C_2\right).
\end{align*}
Thus,
\begin{eqnarray}
{SURE}_i &\ge & 2\sigma^2 [1-\E (Z_i) - s_i\{\E(Z_i)\}^2] \nonumber\\
&\ge &2\sigma^2 \left\{1 -\theta_i(\tilde C_1+ \tilde C_2) \frac{(1+s_i)}{s_i^2}  - \theta_i^2(\tilde C_1+ \tilde C_2)^2\frac{(1+s_i)^2}{s_i^3}  \right\} . \label{eq:alower}
\end{eqnarray}
Thus, combining Equations (\ref{eq:aupper}) and (\ref{eq:alower}) we get 
$$
 \left\{1 -\theta_i(\tilde C_1+ \tilde C_2) \frac{(1+s_i)}{s_i^2}  - \theta_i^2(\tilde C_1+ \tilde C_2)^2\frac{(1+s_i)^2}{s_i^3}  \right\} \le\frac{{SURE}_i}{2\sigma^2} \le \left\{1+ 2 \theta_i (1+s_i) \left( \frac{C_1}{s_i^2} + \frac{C_2}{s_i^{3/2}}\right) \right\},
$$
for $s_i\ge1, \theta_i\ge 1$. 

\subsection{Proof of Theorem~\ref{th:gl}}\label{app:thgl}
Our proof is similar to the proof of Theorem 1 of \citet{polson2010shrink}. Note from Equations (\ref{eq:alphahat1}--\ref{eq:alpha1}) that integrating out $\alpha_i$ we have
$$
\hat \alpha_i \mid \lambda_i^2, \sigma^2, \tau^2 \stackrel{ind}\sim \Nor(0, \sigma^2(d_i^{-2} + \tau^2\lambda_i^2)).
$$
Let $p(\lambda_i^2) \sim (\lambda_i^2)^{a-1} L(\lambda_i^2)$, as $\lambda_i^2\to \infty$ where $a\le 0$. Define $u_i = \sigma^2(d_i^{-2} + \tau^2\lambda_i^2)$. Then, as in Theorem 1 of \citet{polson2010shrink}, we have
$$
p(u_i) \sim u_i^{a-1} L(u_i),\; \text{as}\; u_i\to\infty.
$$
The marginal of $\hat\alpha_i$ is then given by
$$
m(\hat\alpha_i) = \int \frac{1}{\sqrt{2\pi u_i}} \exp\{-\hat\alpha_i^2/(2 u_i)\} p(u_i) d u_i.
$$
An application of  Theorem 6.1 of \citet{barndorff1982normal} shows that
$$
m(\hat \alpha_i) \sim |\hat\alpha_i|^{2a-1} L(|\hat\alpha_i|)\; \text{as}\; |\hat\alpha_i|\to\infty.
$$
Thus, for large $|\hat\alpha_i|$
\begin{eqnarray}
\frac{\partial \log m(\hat\alpha_i)}{\partial \hat\alpha_i} = \frac{(2a-1)}{|\hat\alpha_i|} + \frac{\partial \log L(|\hat\alpha_i|)}{\partial \hat\alpha_i}. \label{eq:app}
\end{eqnarray}
Clearly, the first term in Equation (\ref{eq:app}) goes to zero as $|\hat\alpha_i| \to\infty$. For the second term, we need to invoke the celebrated representation theorem by Karamata. A proof can be found in \citet{bingham1989regular}.
\begin{result} (Karamata's representation theorem).
A function $L$ is slowly varying if and only if there exists $B > 0$ such that for all $x \ge B$ the function can be written in the form
$$
 L(x)=\exp \left(\eta (x)+\int _{B}^{x}{\frac  {\varepsilon (t)}{t}}\,dt\right),
$$
where
$\eta(x)$ is a bounded measurable function of a real variable converging to a finite number as $x$ goes to infinity
$\varepsilon(x)$  is a bounded measurable function of a real variable converging to zero as $x$ goes to infinity.
\end{result}
Thus, using the properties of $\eta(x)$ and $\varepsilon(x)$ from the result above
$$
\frac{d\log(L(x))}{dx} =  \eta' (x)+ {\frac  {\varepsilon (x)}{x}} \to 0\quad \text{as}\quad x\to\infty. 
$$
Using this in Equation (\ref{eq:app}) shows ${\partial \log m(\hat\alpha_i)}/{\partial \hat\alpha_i} \to 0$ as $|\hat\alpha_i| \to \infty$. By similar calculations, ${\partial^2 \log m(\hat\alpha_i)}/{\partial^2 \hat\alpha_i} \to 0$ as $|\hat\alpha_i| \to \infty$. From Equation (\ref{eq:sure})
\begin{align*}
SURE_i =&  \sigma^4 d_i^{-2} \left \{ \frac{\partial}{\partial \hat \alpha_i} \log m(\hat \alpha_i)\right \}^2 + 2\sigma^2 \left\{1 + \sigma^2d_i^{-2} \frac{\partial^2}{\partial \hat \alpha_i^2 } \log m(\hat \alpha_i)\right\}.
\end{align*}
Thus, $SURE_i \to 2\sigma^2$, almost surely, as $|\hat\alpha_i| \to \infty$.

\subsection{Proof of Theorem \ref{th:small}}\label{app3}
The proof of Theorem~\ref{th:small} makes use of technical lemmas  in Appendix~\ref{app4}. 

Recall from Appendix \ref{app1} that if we define $Z_i = 1/(1 + \tau^2 \lambda_i^2 d_i^2)$ then the density of $Z_i$ is given by
\begin{eqnarray}
(Z_i \mid \hat\alpha_i, d_i, \tau, \sigma^2) \sim \mathrm{CCH} \left(Z_i \mid   1, \frac{1}{2} , 1, \frac{\hat \alpha_i^2 d_i^2}{2\sigma^2} , 1, \frac{1}{\tau^2 d_i^2}\right).\label{eq:Zi}
\end{eqnarray}
Then SURE is given by $SURE=\sum_{i=1}^{n} SURE_i$ with
\begin{eqnarray}
SURE_i &=& 2\sigma^2 [1-\E (Z_i)+ 2s_i\E (Z_i^2)-s_i \{\E (Z_i)\}^2] \nonumber\\
&=& 2\sigma^2 [1-\E (Z_i) + s_i\E (Z_i^2) + s_i \mathrm{Var}(Z_i)],  \label{eq:appR}
\end{eqnarray}
where $s_i = \hat \alpha_i^2 d_i^2/2 \sigma^2$. Thus,
\begin{eqnarray}
\frac{\partial \{SURE_i\}}{\partial s_i} &=& -2\sigma^2 \frac{\partial \E(Z_i)}{\partial s_i} +2\sigma^2 \frac{\partial}{\partial s_i} \{s_i \E(Z_i^2)\} +2\sigma^2 \frac{\partial}{\partial s_i} \{s_i \mathrm{Var}(Z_i)\}\nonumber\\
&:=& \mathrm{I} + \mathrm{II} + \mathrm{III}. \label{eq:derivr}
\end{eqnarray}
Now, as a corollary to Lemma~\ref{lemma:derivs}, $({\partial}/{\partial s_i})\E (Z_i) = \{\E(Z_i)\}^2 - \E(Z_i^2) = -\mathrm{Var} (Z_i) <0$, giving $\mathrm{I} >0$. The strict inequality follows from the fact that $Z_i$ is not almost surely a constant for any $s_i \in \mathbb{R}$ and $({\partial}/{\partial s_i})\E (Z_i)$ is continuous at $s_i=0$. Next, consider $\mathrm{II}$. Define $\theta_i = (\tau^2 d_i^2)^{-1}$ and let $0\le s_i \le1$. Then,
\begin{eqnarray*}
\frac{\partial}{\partial s_i} \{s_i \E(Z_i^2)\}&=& \E(Z_i^2) + s_i \frac{\partial}{\partial s_i} \E(Z_i^2)\\
&=& \E(Z_i^2) + s_i \{\E(Z_i)\E(Z_i^2) - \E(Z_i^{3})\} \qquad \text{(by Lemma~\ref{lemma:derivs})} \\
&=& s_i \E(Z_i)\E(Z_i^2) + \{\E(Z_i^2) - s_i \E(Z_i^{3})\}.
\end{eqnarray*}
Now, clearly, the first term, $s_i \E(Z_i)\E(Z_i^2) \ge 0$. We also have $Z_i^2 - s_i Z_i^3 = Z_i^2( 1- s_i Z_i) \ge 0$ a.s. when $0\le Z_i\le 1$ a.s. and $0\le s_i\le 1$. Thus, the second term $\E(Z_i^2) - s_i \E(Z_i^{3}) \ge 0$. Putting the terms together gives $\mathrm{II} \ge 0$. Finally, consider $\mathrm{III}$. Denote $\E(Z_i) = \mu_i$. Then,
\begin{eqnarray*}
	\frac{\partial}{\partial s_i} \{s_i \mathrm{Var}(Z_i)\}&=& \mathrm{Var} (Z_i) + s_i\frac{\partial}{\partial s_i} \{\mathrm{Var}(Z_i)\}\\
	&=& \mathrm{Var} (Z_i) - s_i \frac{\partial^2 \E(Z_i)}{\partial s_i^2} \\
&=& \E\{(Z_i-\mu_i)^2\} - s_i \E\{(Z_i -\mu_i)^3\} \qquad\text{(by Lemma~\ref{lemma:2derivs})}\\
&=& \E[(Z_i - \mu_i)^2\{1 - s_i(Z_i - \mu_i)\}].
\end{eqnarray*}
Now, $(Z_i - \mu_i)^2\{1 - s_i(Z_i - \mu_i)\} \ge 0$ a.s. when $0\le Z_i\le 1$ a.s. and $0\le s_i\le 1$ and thus, $\mathrm{III}\ge 0$. Using $\mathrm{I}, \mathrm{II}$ and $\mathrm{III}$ in Equation (\ref{eq:derivr}) yields $SURE_i$ is an increasing function of $s_i$ when $0\leq s_i \leq 1$, completing the proof of Part A. 

To prove Part B, we need to derive an upper bound on SURE when $s_i=0$. First, consider $s_i=0$ and $0 < \theta_i \le 1$. we have from Equation (\ref{eq:appR}) that $SURE_i = 2\sigma^2 (1-\E Z_i)$. By Lemma \ref{lemma:derivtheta}, $({\partial}/{\partial\theta_i})\E (Z_i) >0$ and $SURE_i$ is a monotone decreasing function of $\theta_i$, where $\theta_i = (\tau^2 d_i^2)^{-1}$. Next consider the case where $s_i=0$ and $\theta_i \in (1, \infty)$. Define $\tilde Z_i = 1-Z_i \in (0,1)$ when $Z_i \in (0,1)$. Then, by Equation (\ref{eq:Zi}) and a formula on Page 9 of \citet{gordy98}, we have that $\tilde Z_i$ also follows a CCH distribution. Specifically, 
\begin{eqnarray*}
(\tilde Z_i \mid \hat\alpha_i, d_i, \tau, \sigma^2) \sim \mathrm{CCH} \left(\tilde Z_i \mid   \frac{1}{2},1 , 1, -\frac{\hat \alpha_i^2 d_i^2}{2\sigma^2} , 1, {\tau^2 d_i^2}\right),
\end{eqnarray*}
and we have $SURE_i=2\sigma^2\E(\tilde Z_i)$. Define $\tilde \theta_i= \theta_i^{-1}=\tau^2d_i^2$. Then by Lemma \ref{lemma:derivtheta}, $({\partial}/{\partial\tilde\theta_i})\E (\tilde Z_i) =-\mathrm{Cov} (\tilde Z_i,\tilde W_i)>0$ on $0<\tilde\theta_i < 1$. Therefore, $SURE_i$ is a monotone increasing function of $\tilde{\theta_i}$ on $0<\tilde{\theta_i}<1$, or equivalently a monotone decreasing function of $\theta_i$ on $\theta_i \in (1, \infty)$. 

Thus, combining the two cases above, we get that SURE at $s_i=0$ is a monotone decreasing function of $\theta_i$ for any $\theta_i \in (0, \infty)$, or equivalently, an increasing function of $\tau^2 d_i^2$.  Since $0\le \tilde Z_i \le 1$ almost surely, a natural upper bound on $SURE_i$ is $2\sigma^2$. However, it is possible to do better provided $\tau$ is chosen sufficiently small. Assume that $\tau^2 \le d_i^{-2}$.  Then, since $SURE_i$ is monotone increasing in $\theta_i$, the upper bound on SURE is achieved when $\theta_i=(\tau^2 d_i^2)^{-1}=1$. In this case, $\E(Z_i)$ has a particularly simple expression, given by
\begin{align}
\E (Z_i) =& \frac{\int_0^1 z_i(1-z_i)^{-\frac{1}{2}}\{\theta_i+(1-\theta_i)z_i\}^{-1} d z_i}
{\int_0^1\ (1-z_i)^{-\frac{1}{2}}\{\theta_i+(1-\theta_i)z_i\}^{-1} d z_i} \nonumber \\
	=& \frac{\int_0^1 z_i(1-z_i)^{-\frac{1}{2}} d z_i} 
	{\int_0^1\ (1-z_i)^{-\frac{1}{2}} d z_i}=\frac{2}{3}.\label{eq:z}
\end{align}	
Thus, $\sup SURE_i = 2\sigma^2 (1-\E Z_i) = 2\sigma^2/3$
, completing the proof of Part B.

To prove Part C, we first note that when $s_i=1$ we have
$$
SURE_i = 2\sigma^2[1 - \E(Z_i)\vert_{s_i=1} + 2\E(Z_i^2)\vert_{s_i=1}  - \{\E(Z_i)\vert_{s_i=1} \}^2] 
$$
where $E(Z_i)$ and $E(Z_i^2)$ are evaluated at $s_i=1$. Recall that when $\theta_i\ge 1$ and $z_i\in (0,1)$ we have $\theta_i^{-1} \le \{\theta_i + (1-\theta_i)z_i\}^{-1} \le 1$. Thus,
\begin{align}
\E (Z_i^2)\vert_{s_i=1}  =& \frac{\int_0^1 z_i^2(1-z_i)^{-\frac{1}{2}}\{\theta_i+(1-\theta_i)z_i\}^{-1}\exp(- z_i) d z_i}
{\int_0^1 (1-z_i)^{-\frac{1}{2}}\{\theta_i+(1-\theta_i)z_i\}^{-1}\exp(-z_i) d z_i} \nonumber\\
& \le \frac{\int_0^1 z_i^2(1-z_i)^{-\frac{1}{2}}\exp(- z_i) d z_i}
{\theta_i^{-1} \int_0^1 (1-z_i)^{-\frac{1}{2}}\exp(-z_i) d z_i} \approx \theta_i \frac{0.459}{1.076} = 0.43\theta_i, \label{eq:z1}
\end{align}
and
\begin{align}
\E (Z_i)\vert_{s_i=1}  =& \frac{\int_0^1 z_i(1-z_i)^{-\frac{1}{2}}\{\theta_i+(1-\theta_i)z_i\}^{-1}\exp(- z_i) d z_i}
{\int_0^1 (1-z_i)^{-\frac{1}{2}}\{\theta_i+(1-\theta_i)z_i\}^{-1}\exp(-z_i) d z_i},\nonumber \\
& \ge \frac{\theta_i^{-1}  \int_0^1 z_i(1-z_i)^{-\frac{1}{2}}\exp(- z_i) d z_i}
{\int_0^1 (1-z_i)^{-\frac{1}{2}}\exp(-z_i) d z_i} \approx \theta_i^{-1} \frac{0.614}{1.076} = 0.57\theta_i^{-1}. \label{eq:z21}
\end{align}
Thus,
$$
SURE_i \le 2\sigma^2\left[1 - \frac{0.57}{\theta_i} + 0.86\theta_i  - \left(\frac{0.57}{\theta_i}\right)^2\right]. 
$$
When $\theta_i=1$, it can be seen that $SURE_i\le 1.93\sigma^2$.

\subsection{Proof of Theorem \ref{th:riskgl}}\label{app:riskgl}
The proof of Theorem~\ref{th:riskgl} makes use of technical lemmas  in Appendix~\ref{app4}. 

Recall from Appendix \ref{app1} that if we define $Z_i = 1/(1 + \tau^2 \lambda_i^2 d_i^2)$ then the density of $Z_i$ is given by
\begin{eqnarray}
(Z_i \mid \hat\alpha_i, d_i, \tau, \sigma^2) \sim \mathrm{CCH} \left(Z_i \mid   1, {1}/{2} , 1, s_i, 1, \theta_i \right).\label{eq:Zi}
\end{eqnarray}
where $s_i = \hat \alpha_i^2 d_i^2/2 \sigma^2$ and $\theta_i = (\tau^2 d_i^2)^{-1}$. Consider the case where $d_i = 1$ for all $i$ and $\tau^2 = 1$, i.e., $\theta_i=1$ for all $i$. From Equation (\ref{eq:appR}), the risk estimate is $SURE=\sum_{i=1}^{n} SURE_i$ with
\begin{eqnarray*}
SURE_i &=& 2\sigma^2 [1-\E (Z_i) + s_i\E (Z_i^2) + s_i \mathrm{Var}(Z_i)], \\
&\le&  2\sigma^2 [1-\E (Z_i) + s_i + s_i \mathrm{Var}(Z_i)] = \check{R}_i.
\end{eqnarray*}
We begin by showing that the upper bound $ \check{R}_i = 2\sigma^2 [1-\E (Z_i) + s_i + s_i \mathrm{Var}(Z_i)]$ is convex in $s_i$ when $s_i\in (0,1)$. It suffices to show $-\E(Z_i)$ and $s_i \mathrm{Var}(Z_i)$ are separately convex. First, $({\partial^2}/{\partial^2 s_i})\E (Z_i) =\E\{(Z_i-\mu_i)^3\} \le 0$, by Lemmas \ref{lemma:2derivs} and \ref{lemma:skew}, proving $-\E(Z_i)$ is convex. Next,
\begin{eqnarray*}
	\frac{\partial^2}{\partial s_i^2} \{s_i \mathrm{Var}(Z_i)\}&=& \frac{\partial}{\partial s_i} \left[\mathrm{Var} (Z_i) + s_i\frac{\partial}{\partial s_i} \{\mathrm{Var}(Z_i)\}\right]\\
	&=& 2 \frac{\partial}{\partial s_i} \{\mathrm{Var}(Z_i)\} + s_i \frac{\partial^2}{\partial s_i^2} \{\mathrm{Var}(Z_i)\}\\
	&=& -2 \E(Z_i - \mu_i)^3 - s_i \frac{\partial}{\partial s_i} \E(Z_i - \mu_i)^3 \quad \text{(by Lemma~\ref{lemma:2derivs})}\\
	&=& -2 \E(Z_i - \mu_i)^3 + s_i \E(Z_i - \mu_i)^4, \quad \text{(by Lemma~\ref{lemma:3derivs})}\\
	&\ge& 0,
\end{eqnarray*}
where the last inequality follows by Lemma~\ref{lemma:skew}. Thus, since $\check{R}_i$ is convex, it lies entirely below the straight line joining the two end points for $s_i \in (0,1)$. But $\check{R}_i \vert_{s_i=0} \le 2\sigma^2/3 = 0.67\sigma^2$ (by Equation (\ref{eq:z})) and 
$$
\check{R}_i \vert_{s_i=1} \le 2\sigma^2\left[1 - {0.57} + 1 + 0.43  - (0.57)^2\right] = 3.07\sigma^2,
$$
by Equations (\ref{eq:z1}) and (\ref{eq:z21}).Thus, by convexity 
\begin{equation}
SURE_i \le \check{R}_i \le 0.67\sigma^2 + s_i (3.07 - 0.67)\sigma^2 = (0.67+ 2.4s_i)\sigma^2 \text{ for } s_i \in (0,1) \label{eq:bound01}
\end{equation}
We remark here that our simulations suggest $SURE_i$ itself is convex, not just the upper bound $\check{R}_i$, although a proof seems elusive. Nevertheless, as we shall see below, the convexity of $\check{R}_i$ is sufficient for our purposes.

Next, consider the interval $s_i \in (1,3)$. Noting that both $\E(Z_i)$ and $\E(Z_i^2)$ are monotone decreasing functions of $s_i$ we have 
$$
SURE_i \le 2\sigma^2[1 - \E(Z_i)\vert_{s_i=3} + 2s_i\{\E(Z_i^2)\vert_{s_i=1} \} - s_i\{\E(Z_i)\vert_{s_i=3} \}^2] 
$$
But,
\begin{align}
\E (Z_i)\vert_{s_i=3, \theta_i = 1}  =& \frac{\int_0^1 z_i(1-z_i)^{-\frac{1}{2}}\exp(- 3z_i) d z_i}
{\int_0^1 (1-z_i)^{-\frac{1}{2}}\exp(-3z_i) d z_i} = 0.35. \nonumber
\end{align}
$\E (Z_i^2)\vert_{s_i=1}<0.43$ from Equation (\ref{eq:z1}). Thus,
\begin{equation}
SURE_i \le 2\sigma^2[1 - 0.35 + 0.86s_i - s_i(0.35)^2 \}^2] = 2\sigma^2(0.65 + 0.74s_i) \text{  for  } s_i \in (1,3).\label{eq:bound13}
\end{equation}
Using the upper bound from Theorem~\ref{th:large}, 
\begin{equation}
SURE_i \le 11.55\sigma^2 \text{ for } s_i \ge 3\label{eq:bound3}.
\end{equation}
When $\alpha_i = 0$, we have that $\hat\alpha_i \sim \Nor(0, \sigma^2 d_i^{-2})$. Thus, $\hat\alpha_i^2 d_i^2/\sigma^2 \sim \chi^2(1)$. Since $s_i =\hat \alpha_i^2 d_i^2/2\sigma^2$ we have that $p(s_i) = (\pi)^{-1/2}s_i^{-1/2} \exp(-s_i)$ for $s_i \in (0,\infty)$. Combining Equations (\ref{eq:bound01}), (\ref{eq:bound13}) and (\ref{eq:bound3}) we have
\begin{align*}
\mathrm{Risk}_i = \E(SURE_i) &\le \int_{0}^{1} \sigma^2(0.67+ 2.4s_i) \pi^{-1/2} s_i^{-1/2} \exp(-s_i) ds_i\\
&+ \int_{1}^{3} 2\sigma^2(0.65+ 0.74s_i) \pi^{-1/2} s_i^{-1/2} \exp(-s_i) ds_i\\
&+ \int_{3}^{\infty} 11.55\sigma^2 \pi^{-1/2} s_i^{-1/2} \exp(-s_i) ds_i\\
&= 1.75\sigma^2.
\end{align*}

\subsection{Technical lemmas}\label{app4}
\begin{lemma} \label{lemma:derivs}
	If $Z \sim \mathrm{CCH}(p,q,r,s, \nu ,\theta)$, then $({\partial}/{\partial s})\E (Z^k)= \E(Z)\E(Z^k) - \E(Z^{k+1})$.
\end{lemma}
\begin{lemma} \label{lemma:2derivs}
	If $Z \sim \mathrm{CCH}(p,q,r,s,\nu ,\theta)$, then $({\partial^2}/{\partial^2 s})\E (Z)= -({\partial}/{\partial s})\mathrm{Var} (Z)= \E\{(Z-\mu)^3\}$, where $\mu = \E(Z)$.
\end{lemma}
\begin{lemma} \label{lemma:derivtheta}
	If $Z \sim \mathrm{CCH}(p,q,r,s,\nu,\theta)$, then $({\partial}/{\partial \theta})\E (Z) =-\mathrm{Cov} (Z,W)$, for $W=(1-\nu Z)\{\theta+(1-\theta)\nu Z\}^{-1}$. If $0<\theta\le1$ then $({\partial}/{\partial \theta})\E (Z)>0$.
\end{lemma}
\begin{lemma} \label{lemma:skew}
	If $Z \sim \mathrm{CCH}(p,q,r,s, 1 ,1)$ with $q>p$, then $\E(Z - \mu)^3 \le 0$, where $\mu = \E(Z)$.
\end{lemma}
\begin{lemma} \label{lemma:3derivs}
	If $Z \sim \mathrm{CCH}(p,q,r,s,\nu ,\theta)$, then $({\partial}/{\partial s})\E (Z - \mu)^3= -\E\{(Z-\mu)^4\}$, where $\mu = \E(Z)$.
\end{lemma}

\subsubsection{Proof of Lemma \ref{lemma:derivs}}
Let, $Z \sim \mathrm{CCH} (p,q,r,s,\nu,\theta)$. Then for any integer $k$
\begin{align*}
	\E (Z^k) =& \frac{\int_0^{1/\nu} z^{k+p-1}(1- \nu z)^{q-1}\{\theta+(1-\theta) \nu z\}^{-r}\exp(-sz) dz}
	{\int_0^{1/\nu} z^{p-1}(1- \nu z)^{q-1}\{\theta+(1-\theta) \nu z\}^{-r}\exp(-sz) dz}.
	\end{align*}
	Thus,
	\begin{align*}
	\frac{\partial}{\partial s}\E (Z^k) = & 
	\frac{\int_0^{1/\nu} -z^{k+p}(1-\nu z)^{q-1}\{\theta+(1-\theta)\nu z\}^{-r}\exp(-sz) dz}
	{\int_0^{1/\nu} z^{p-1}(1-\nu z)^{q-1}\{\theta+(1-\theta) \nu z\}^{-r}\exp(-sz) dz} \\
	-& \Bigg[\frac{\int_0^{1/\nu} z^{k+p-1}(1-\nu z)^{q-1}\{\theta+(1-\theta) \nu z\}^{-r}\exp(-sz) dz}
	{\int_0^{1/\nu} z^{p-1}(1-\nu z)^{q-1}\{\theta+(1-\theta)\nu z\}^{-r}\exp(-sz) dz}\\
	&\times \frac{\int_0^{1/\nu} -z^{p}(1-\nu z)^{q-1}\{\theta+(1-\theta) \nu z\}^{-r}\exp(-sz) dz}
	{\int_0^{1/\nu} z^{p-1}(1-\nu z)^{q-1}\{\theta+(1-\theta) \nu z\}^{-r}\exp(-sz) dz}\Bigg]	\\
	=& -\E (Z^{k+1})+ \E(Z)\E(Z^k).
\end{align*}
For an alternative proof directly using the $H(\cdot)$ functions, see Appendix D of \citet{gordy98}.

\subsubsection{Proof of Lemma \ref{lemma:2derivs}}
Let, $Z \sim \mathrm{CCH} (p,q,r,s,\nu,\theta)$. From Lemma~\ref{lemma:derivs}, $({\partial}/{\partial s})\E (Z) = -\E(Z^2)  + \{\E(Z)\}^2 = -\mathrm{Var}(Z)$. Let $\mu = \E(Z)$. Then,
\begin{align*}
\frac{\partial^2}{\partial s^2}\E (Z) =& -\frac{\partial}{\partial s} \mathrm{Var} (Z)\\
=& -\frac{\partial}{\partial s} \left[\frac{\int_0^{1/\nu} (z - \mu)^2 z^{p-1}(1-\nu z)^{q-1}\{\theta+(1-\theta) \nu z\}^{-r}\exp(-sz) dz}
	{\int_0^{1/\nu} z^{p-1}(1-\nu z)^{q-1}\{\theta+(1-\theta) \nu z\}^{-r}\exp(-sz) dz}\right]\\
=& \frac{\int_0^{1/\nu} (z - \mu)^2 z^{p}(1-\nu z)^{q-1}\{\theta+(1-\theta) \nu z\}^{-r}\exp(-sz) dz}
	{\int_0^{1/\nu} z^{p-1}(1- \nu z)^{q-1}\{\theta+(1-\theta) \nu z\}^{-r}\exp(-sz) dz} \\
	&- \Bigg[\frac{\int_0^{1/\nu} (z-\mu)^2 z^{p-1}(1-\nu z)^{q-1}\{\theta+(1-\theta) \nu z\}^{-r}\exp(-sz) dz}
	{\int_0^{1/\nu} z^{p-1}(1- \nu z)^{q-1}\{\theta+(1-\theta)\nu z\}^{-r}\exp(-sz) dz}\\
	&\quad\times \frac{\int_0^{1/\nu} z^{p}(1- \nu z)^{q-1}\{\theta+(1-\theta) \nu z\}^{-r}\exp(-sz) dz}
	{\int_0^{1/\nu} z^{p-1}(1- \nu z)^{q-1}\{\theta+(1-\theta) \nu z\}^{-r}\exp(-sz) dz}\Bigg]	\\
	=& \mathrm{Cov} (Z, (Z-\mu)^2)\\
		=& \E[(Z-\mu) \{(Z-\mu)^2 - \E(Z-\mu)^2\}]\\
		=& \E\{(Z-\mu)^3\} - \mathrm{Var}(Z)\E(Z-\mu)= \E\{(Z-\mu)^3\} .
\end{align*}

\subsubsection{Proof of Lemma \ref{lemma:derivtheta}}
Let $Z \sim \mathrm{CCH} (p,q,r,s,\nu,\theta)$ and $W=(1-\nu Z)\{\theta+(1-\theta)\nu Z\}^{-1}$. Then,
\begin{align*}
	\frac{\partial}{\partial \theta}\E (Z) = & 
	- \frac{\int_0^{1/\nu} z^{p}(1- \nu z)^{q}\{\theta+(1-\theta) \nu z\}^{-(r+1)}\exp(-sz) dz}
	{\int_0^{1/\nu} z^{p-1}(1-\nu z)^{q-1}\{\theta+(1-\theta)\nu z\}^{-r}\exp(-sz) dz}\\
	 &+ \Bigg[\frac{\int_0^{1/\nu} z^{p}(1- \nu z)^{q-1}\{\theta+(1-\theta) \nu z\}^{-r}\exp(-sz) dz}
	{\int_0^{1/\nu} z^{p-1}(1-\nu z)^{q-1}\{\theta+(1-\theta) \nu z\}^{-r}\exp(-sz) dz} \\
	&\quad \times \frac{\int_0^{1/\nu} z^{p-1}(1- \nu z)^{q}\{\theta+(1-\theta) \nu z\}^{-(r+1)}\exp(-sz) dz}
	{\int_0^{1/\nu} z^{p-1}(1- \nu z)^{q-1}\{\theta+(1-\theta) \nu z\}^{-r}\exp(-sz) dz}	 \Bigg] \\
	=& -\E (ZW) + \E (Z) \E (W) = -\mathrm{Cov}(Z,W).
\end{align*}
When $0<\theta\le1$, it is obvious that $Z$ and $W$ are negatively correlated, and thus $-\mathrm{Cov}(Z,W)>0$. 

\subsubsection{Proof of Lemma \ref{lemma:skew}}
Let $Z \sim \mathrm{CCH} (p,q,r,s,1, 1)$. Then,
\begin{align*}
\E (Z - \mu)^3  =& \frac{\int_0^{1} (z - \mu)^3 z^{p-1}(1- z)^{q-1}\exp(-sz) dz}
	{\int_0^{1} z^{p-1}(1- z)^{q-1}\exp(-sz) dz},
\end{align*}
which can be seen to have the same sign as the third central moment, or skewness of a $\mathrm{Beta}(p,q)$ random variable, which is negative when $q>p$.

\subsubsection{Proof of Lemma \ref{lemma:3derivs}}

Let, $Z \sim \mathrm{CCH} (p,q,r,s,\nu,\theta)$. Let $\mu = \E(Z)$. Then,
\begin{align*}
\frac{\partial}{\partial s} \E(Z-\mu)^3   =& -\frac{\int_0^{1/\nu} (z - \mu)^3 z^{p}(1-\nu z)^{q-1}\{\theta+(1-\theta) \nu z\}^{-r}\exp(-sz) dz}
	{\int_0^{1/\nu} z^{p-1}(1- \nu z)^{q-1}\{\theta+(1-\theta) \nu z\}^{-r}\exp(-sz) dz} \\
	&+ \Bigg[\frac{\int_0^{1/\nu} (z-\mu)^3 z^{p-1}(1-\nu z)^{q-1}\{\theta+(1-\theta) \nu z\}^{-r}\exp(-sz) dz}
	{\int_0^{1/\nu} z^{p-1}(1- \nu z)^{q-1}\{\theta+(1-\theta)\nu z\}^{-r}\exp(-sz) dz}\\
	&\quad\times \frac{\int_0^{1/\nu} z^{p}(1- \nu z)^{q-1}\{\theta+(1-\theta) \nu z\}^{-r}\exp(-sz) dz}
	{\int_0^{1/\nu} z^{p-1}(1- \nu z)^{q-1}\{\theta+(1-\theta) \nu z\}^{-r}\exp(-sz) dz}\Bigg]	\\
	=& -\mathrm{Cov} (Z, (Z-\mu)^3)\\
		=& -\E[(Z-\mu) \{(Z-\mu)^3 - \E(Z-\mu)^3\}]\\
		=& -\E\{(Z-\mu)^4\} + \E(Z-\mu)^3\E(Z-\mu)= -\E\{(Z-\mu)^4\} .
\end{align*}

\bibliography{horseshoe-plus}

\clearpage\pagebreak\newpage
\begin{center}
{\LARGE{\bf Supplementary Material to\\ {\it Prediction risk for the horseshoe regression}}}
\end{center}
\vskip 2cm
\baselineskip=15pt
\begin{center}
\vspace{-1cm}
Anindya Bhadra\\
Department of Statistics, Purdue University, 250 N. University Street, West Lafayette, IN 47907, USA.\\
bhadra@purdue.edu\\
\hskip 5mm \\
Jyotishka Datta\\
Department of Mathematical Sciences, University of Arkansas, Fayetteville, AR 72701, USA.\\
jd033@uark.edu\\
\hskip 5mm \\
Yunfan Li\\
Department of Statistics, Purdue University, 250 N. University Street, West Lafayette, IN 47907, USA.\\
li896@purdue.edu\\
\hskip 5mm \\
Nicholas G. Polson and Brandon Willard\\
The University of Chicago Booth School of Business, 5807 S. Woodlawn Ave., Chicago, IL 60637, USA.\\
ngp@chicagobooth.edu, brandonwillard@gmail.com\\
\end{center}

\setcounter{equation}{0}
\setcounter{table}{0}
\setcounter{section}{0}
\setcounter{subsection}{0}
\setcounter{figure}{0}
\renewcommand{\theequation}{S.\arabic{equation}}
\renewcommand{\thesection}{S.\arabic{section}}
\renewcommand{\thesubsection}{S.\arabic{subsection}}
\renewcommand{\thetable}{S.\arabic{table}}
\renewcommand{\thefigure}{S.\arabic{figure}}

\newpage
\subsection{Additional simulations} \label{sec:suppsim}
We provide additional simulation results, complementing the results in Table~\ref{tab:result}. For each simulation setting, we report SURE when a formula is available. We also report the average out of sample prediction SSE (standard deviation of SSE) computed based on one training set and 200 testing sets. For each setting, $n=100$. The methods under consideration are ridge regression (RR), principal components regression (PCR), the lasso, the adaptive lasso (A\_LASSO), the minimax concave penalty (MCP), the smoothly clipped absolute deviation (SCAD) penalty and the proposed horseshoe regression (HS). The method with the \textbf{lowest SSE} is in bold and that with \emph{lowest SURE} is in italics for each setting. The features of these additional simulations include the following.
\begin{enumerate}
\item We explore a higher dimensional case ($p=1000$) for each setting.
\item We incorporate two non-convex regression methods for comparisons. These are SCAD \citep{fan2001variable} and MCP \citep{zhang2010nearly}. 
\item We explore different choices of the design matrix $X$. These include three cases: (i) $X$ is generated from a factor model, where it is relatively ill-conditioned (as in Table~\ref{tab:result}), (ii) $X$ is generated from a standard normal, where it is well-conditioned and (iii) $X$ is exactly orthogonal, with all singular values equal to 1. These are reported in corresponding table captions. 
\item We explore different choices of true $\alpha$. These include three cases: (i) Sparse-robust $\alpha$, where most elements of $\alpha$ are close to zero and a few are large, (ii) null $\alpha$, where all elements of $\alpha$ are zero and (iii) dense $\alpha$, where all elements are non-zero. Exact settings and the value of $||\alpha||^2$ are reported in the table captions.
\end{enumerate}

The major finding is that the horseshoe regression outperforms the other global shrinkage methods (ridge and PCR) when $\alpha$ is sparse-robust, which is consistent with the theoretical observation in Section~\ref{sec:risk}.  It also outperforms the other selection-based methods in this case.  On the other hand, the dense $\alpha$ case is most often favorable to ridge regression, while the null $\alpha$ case is favorable to selection-based methods such as the lasso, adaptive lasso, MCP or SCAD, due to the ability of these methods to produce exact zero estimates. However, the selection-based methods perform considerably worse compared to both global and global-local shrinkage methods in the dense $\alpha$ case.

\newpage
\setlength{\oddsidemargin}{-12pt}
\setlength{\evensidemargin}{-12pt}

 \begin{table}[!h]
 \vspace{1cm}
    	\centering
    	\caption{Sparse-robust $\alpha$ (five large coefficients equal to $10$ and other coefficients equal to $0.5$ or $-0.5$ randomly, giving $\sum_{i=1}^{n}{\alpha_i^2} = 523.75$); $X$ generated by a factor model with 4 factors, each factor follows a standard normal distribution; $d_1/d_n$ is the ratio of largest and smallest singular values of $X$.}
    	\label{tab:s1}
	\footnotesize
	\begin{tabular}{c|c|cc|cc|cc|c|c|c|cc}
    	\toprule
		& &\multicolumn{2}{c|}{RR} &\multicolumn{2}{c|}{PCR} & \multicolumn{2}{c|}{LASSO} & A\_LASSO & MCP & SCAD & \multicolumn{2}{c}{HS}\\
    	$p$ & $d_1/d_n$ & SURE & SSE & SURE & SSE & SURE & SSE & SSE & SSE & SSE & SURE & SSE \\
    	\toprule
    	100 & 2360.43 & 165.45 & 159.83 & 163.80 & 161.62 & 122.78 & 145.07 & 132.25 & 127.07 & 127.85 & \textit{116.01} & \textbf{123.07} \\
    	& & & (22.02) & & (21.28) & & (19.39) & (17.57) & (16.71) & (17.19) & & (16.43) \\
    	200 & 28.47 & 188.13 & 206.39 & 217.40 & 244.71 & 174.48 & 162.94 & \textbf{148.41} & 154.01 & 157.73 & \textit{160.89} & 152.37 \\
    	& & & (28.61) & & (29.80) & & (24.44) & (22.48) & (23.17) & (23.23) & & (22.75) \\
    	300 & 22.76 & 192.35 & 212.05 & 266.84 & 280.25 & \textit{155.26} & 190.09 & 175.46 & 172.18 & 176.29 & 157.50 & \textbf{164.17} \\
    	& & & (28.50) & & (32.62) & & (26.20) & (22.17) & (21.55) & (22.19) & & (22.85) \\
    	400 & 21.81 & 194.73 & 199.36 & 337.32 & 328.48 & 179.45 & 182.89 & 197.25 & 199.08 & 198.40 & \textit{172.63} & \textbf{165.15} \\
    	& & & (28.75) & & (34.79) & & (27.41) & (25.02) & (25.52) & (25.31) & & (24.67) \\
    	500 & 18.18 & 196.03 & 180.12 & 410.82 & 379.03 & 158.07 & 173.82 & 223.21 & 224.91 & 226.76 & \textit{166.10} & \textbf{161.77} \\
    	& & & (27.16) & & (39.41) & & (26.50) & (27.98) & (29.65) & (29.26) & & (24.22) \\
    	1000 & 15.20 & 197.91 & 184.86 & 669.69 & 736.69 & 196.83 & 205.28 & 345.26 & 344.04 & 344.04 & \textit{191.64} & \textbf{182.18} \\
    	& & & (26.42) & & (56.58) & & (29.56) & (36.60) & (37.34) & (37.34) & & (25.43) \\
    	\bottomrule
    	\end{tabular}
 \end{table} 

 \begin{table}[!h]
  \vspace{2cm}
    	\centering
    	\caption{Null $\alpha$ ($\sum_{i=1}^{n} \alpha_i^2 = 0$); $X$ is the same as in Table~\ref{tab:s1}.}
    	\label{tab:s2}
\footnotesize
	\begin{tabular}{c|cc|cc|cc|c|c|c|cc}
    	\toprule
		&\multicolumn{2}{c|}{RR} &\multicolumn{2}{c|}{PCR} & \multicolumn{2}{c|}{LASSO} & A\_LASSO & MCP & SCAD & \multicolumn{2}{c}{HS}\\
    	$p$ & SURE & SSE & SURE & SSE & SURE & SSE & SSE & SSE & SSE & SURE & SSE \\
    	\toprule
		100 & \textit{88.23} & 100.86 & 92.85 & 113.28 & 87.36 & 100.81 & \textbf{100.70} & 100.81 & 100.81 & 92.42 & 102.31 \\
		& & (13.20) & & (14.91) & & (13.29) & (13.21) & (13.29) & (13.29) & & (13.72) \\
		200 & 121.30 & 107.68 & 128.83 & 115.65 & \textit{117.90} & 105.77 & \textbf{100.32} & 104.39 & 101.78 & 122.29 & 111.39 \\
		& & (15.70) & & (16.28) & & (15.06) & (14.80) & (14.93) & (14.89) & & (16.12) \\
		300 & 125.78 & 101.36 & 139.96 & 124.37 & \textit{108.85} & 111.85 & \textbf{101.30} & 104.89 & 102.91 & 119.67 & 112.00 \\
		& & (13.99) & & (17.35) & & (15.37) & (14.02) & (14.27) & (14.00) & & (15.76) \\
		400 & 113.00 & 99.50 & 113.50 & 99.41 & \textit{102.81} & 111.92 & 114.62 & \textbf{99.40} & 110.30 & 113.42 & 107.20 \\
		& & (13.12) & & (13.09) & & (15.51) & (15.80) & (13.20) & (15.20) & & (14.90) \\
		500 & 90.74 & 101.04 & \textit{88.26} & 107.31 & 90.26 & 99.49 & \textbf{99.06} & 99.49 & 99.49 & 101.55 & 102.93 \\
		& & (14.17) & & (15.08) & & (14.16) & (14.04) & (14.16) & (14.16) & & (14.68) \\
		1000 & 88.86 & 100.34 & 85.67 & 103.47 & \textit{82.51} & 100.43 & \textbf{99.52} & 100.43 & 100.41 & 99.73 & 104.84 \\
		& & (14.00) & & (14.29) & & (13.90) & (13.70) & (13.90) & (14.00) & & (14.96) \\		
    	\bottomrule
    	\end{tabular}
 \end{table}

 \begin{table}[!h]
    	\centering
    	\caption{Dense $\alpha$ (all coefficients equal to $2$, giving $\sum_{i=1}^{n}{\alpha_i^2} = 400$);  $X$ is the same as in Table~\ref{tab:s1}.}
    	\label{tab:s3}
\footnotesize
	\begin{tabular}{c|cc|cc|cc|c|c|c|cc}
    	\toprule
		&\multicolumn{2}{c|}{RR} &\multicolumn{2}{c|}{PCR} & \multicolumn{2}{c|}{LASSO} & A\_LASSO & MCP & SCAD & \multicolumn{2}{c}{HS}\\
    	$p$ & SURE & SSE & SURE & SSE & SURE & SSE & SSE & SSE & SSE & SURE & SSE \\
    	\toprule
    	100 & \textit{162.49} & \textbf{159.94} & 177.47 & 175.19 & 194.86 & 203.89 & 504.55 & 491.67 & 491.67 & 185.46 & 173.89 \\
    	& & (21.60) & & (22.36) & & (28.11) & (46.13) & (45.31) & (45.31) & & (23.63) \\
    	200 & \textit{183.75} & \textbf{200.92} & 196.06 & 233.12 & 211.99 & 232.36 & 960.77 & 895.83 & 911.94 & 204.10 & 228.18 \\
    	& & (27.97) & & (31.18) & & (31.06) & (59.48) & (60.85) & (60.19) & & (31.21) \\
    	300 & \textit{189.38} & \textbf{209.92} & 200.39 & 225.92 & 216.01 & 524.84 & 1344.27 & 1298.80 & 1298.80 & 206.99 & 227.55 \\
    	& & (27.88) & & (30.15) & & (69.45) & (71.29) & (77.97) & (77.97) & & (29.98) \\
    	400 & \textit{193.02} & \textbf{195.01} & 197.74 & 217.68 & 218.16 & 306.15 & 1768.05 & 1675.73 & 1675.73 & 207.81 & 213.91 \\
    	& & (28.92) & & (31.02) & & (42.65) & (78.92) & (75.86) & (75.86) & & (31.14) \\
    	500 & \textit{194.85} & \textbf{175.46} & 208.87 & 201.18 & 220.34 & 743.40 & 2154.61 & 2082.54 & 2081.42 & 207.93 & 188.93 \\ 
    	& & (26.52) & & (29.07) & & (100.54) & (92.70) & (92.93) & (93.37) & & (28.10) \\
    	1000 & \textit{197.37} & \textbf{181.65} & 247.40 & 197.59 & 224.75 & 210.78 & 4280.80 & 4075.00 & 4075.00 & 203.47 & 186.48 \\
    	& & (26.50) & & (27.76) & & (29.70) & (145.28) & (138.72) & (138.72) & & (26.95) \\
    	\bottomrule
    	\end{tabular}
 \end{table}

 \begin{table}[!h]
    	\centering
    	\caption{Sparse-robust $\alpha$ (five large coefficients equal to $10$ and other coefficients equal to $0.5$ or $-0.5$ randomly, giving $\sum_{i=1}^{n}{\alpha_i^2} = 523.75$); $X$ follows a standard normal distribution; $d_1/d_n$ is the ratio of largest and smallest singular values of $X$.}
    	\label{tab:s4}
\footnotesize
	\begin{tabular}{c|c|cc|cc|cc|c|c|c|cc}
    	\toprule
		& & \multicolumn{2}{c|}{RR} &\multicolumn{2}{c|}{PCR} & \multicolumn{2}{c|}{LASSO} & A\_LASSO & MCP & SCAD & \multicolumn{2}{c}{HS}\\
    	$p$ & $d_1/d_n$ & SURE & SSE & SURE & SSE & SURE & SSE & SSE & SSE & SSE & SURE & SSE \\
    	\toprule
    	100 & 351.2 & 196.72 & 188.67 & 228.78 & 231.34 & 207.63 & 425.67 & 2537.23 & 2573.27 & 2573.27 & \textit{195.22} & \textbf{188.52} \\
    	& & & (29.04) & & (34.36) & & (59.68) & (112.58) & (128.46) & (128.46) & & (28.90) \\
    	200 & 5.73 & \textit{199.84} & \textbf{193.41} & 221.35 & 206.25 & 218.26 & 1618.40 & 4849.94 & 4915.72 & 4964.26 & 201.91 & 194.14 \\
    	& & & (28.36) & & (28.28) & & (211.54) & (146.61) & (186.74) & (186.55) & & (28.45) \\
    	300 & 3.63 & \textit{199.91} & \textbf{217.43} & 8538.46 & 8082.93 & 222.62 & 1926.13 & 13132.01 & 7316.38 & 7316.39 & 200.92 & 219.97 \\
    	& & & (27.91) & & (320.98) & & (248.97) & (281.12) & (218.82) & (218.83) & & (28.13) \\
    	400 & 2.89 & \textit{199.94} & \textbf{197.47} & 228.38 & 223.43 & 224.53 & 2384.04 & 9593.41 & 9695.47 & 9695.47 & 200.31 & 197.69 \\
    	& & & (27.43) & & (31.13) & & (299.58) & (210.42) & (323.72) & (323.72) & & (27.46) \\
    	500 & 2.51 & \textit{199.95} & \textbf{193.86} & 256.09 & 273.63 & 224.96 & 471.73 & 11980.15 & 11991.11 & 11991.11 & 200.15 & 194.17 \\
    	& & & (27.26) & & (33.97) & & (60.52) & (235.77) & (272.80) & (272.80) & & (27.24) \\
    	1000 & 1.88 & 199.98 & 185.85 & 605.64 & 560.45 & 222.18 & 5781.13 & 23866.06 & 24566.13 & 24566.13 & \textit{199.96} & \textbf{185.79} \\
    	& & & (27.17) & & (45.94) & & (759.08) & (326.06) & (941.16) & (941.16) & & (27.17) \\
    	\bottomrule
    	\end{tabular}
 \end{table} 
 
  \begin{table}[!h]
    	\centering
    	\caption{Null $\alpha$ ($\sum_{i=1}^{n} \alpha_i^2 = 0$); $X$ is the same as in Table~\ref{tab:s4}.}
    	\label{tab:s5}
\footnotesize
	\begin{tabular}{c|cc|cc|cc|c|c|c|cc}
    	\toprule
		&\multicolumn{2}{c|}{RR} &\multicolumn{2}{c|}{PCR} & \multicolumn{2}{c|}{LASSO} & A\_LASSO & MCP & SCAD & \multicolumn{2}{c}{HS}\\
    	$p$ & SURE & SSE & SURE & SSE & SURE & SSE & SSE & SSE & SSE & SURE & SSE \\
    	\toprule
    	100 & 118.45 & 119.12 & 96.35 & 106.88 & \textit{92.06} & 101.21 & \textbf{100.52} & 101.21 & 101.21 & 119.11 & 114.69 \\
    	& & (18.19) & & (15.18) & & (14.33) & (14.20) & (14.33) & (14.33) & & (17.47) \\
    	200 & 136.93 & 135.02 & 96.49 & 100.14 & \textit{94.54} & 100.15 & \textbf{100.13} & 100.39 & 102.06 & 126.19 & 126.34 \\
    	& & (21.74) & & (14.77) & & (14.69) & (14.70) & (14.88) & (15.26) & & (20.13) \\
    	300 & 152.52 & 160.29 & 119.00 & 131.94 & \textit{118.15} & 100.71 & \textbf{100.49} & 100.71 & 100.71 & 140.91 & 140.08 \\
    	& & (21.61) & & (18.01) & & (14.48) & (14.37) & (14.48) & (14.48) & & (18.93) \\
    	400 & 158.64 & 159.13 & 100.88 & 104.06 & \textit{96.30} & 103.07 & \textbf{100.46} & 103.03 & 103.03 & 138.62 & 132.14 \\
    	& & (23.15) & & (15.59) & & (15.11) & (14.82) & (15.04) & (15.04) & & (19.62) \\
    	500 & 166.06 & 158.83 & 98.64 & 98.10 & \textit{94.30} & 100.36 & \textbf{97.99} & 100.36 & 100.36 & 140.14 & 131.53 \\
    	& & (23.35) & & (14.50) & & (14.79) & (14.50) & (14.79) & (14.79) & & (19.59) \\
    	1000 & 181.23 & 169.22 & 89.95 & 100.66 & \textit{87.79} & 100.07 & \textbf{99.80} & 100.66 & 100.51 & 141.11 & 138.94 \\
    	& & (25.25) & & (14.10) & & (14.03) & (14.00) & (14.08) & (14.07) & & (21.12) \\
    	\bottomrule
    	\end{tabular}
 \end{table} 
 
   \begin{table}[!h]
    	\centering
    	\caption{Dense $\alpha$ (all coefficients equal to $2$, giving $\sum_{i=1}^{n}{\alpha_i^2} = 400$); $X$ is the same as in Table~\ref{tab:s4}.}
    	\label{tab:s6}
\footnotesize
	\begin{tabular}{c|cc|cc|cc|c|c|c|cc}
    	\toprule
		&\multicolumn{2}{c|}{RR} &\multicolumn{2}{c|}{PCR} & \multicolumn{2}{c|}{LASSO} & A\_LASSO & MCP & SCAD & \multicolumn{2}{c}{HS}\\
    	$p$ & SURE & SSE & SURE & SSE & SURE & SSE & SSE & SSE & SSE & SURE & SSE \\
    	\toprule
    	100 & \textit{193.13} & \textbf{188.60} & 206.31 & 200.53 & 222.52 & 210.23 & 40019.73 & 40063.42 & 40063.42 & 199.25 & 191.74 \\
	   	& & (28.91) & & (29.51) & & (31.36) & (690.71) & (717.27) & (717.27) & & (29.37) \\
    	200 & \textit{199.76} & \textbf{193.93} & 392.38 & 349.52 & 224.73 & 316.05 & 80016.11 & 80187.49 & 80187.49 & 200.14 & 194.48 \\
    	& & (28.41) & & (37.90) & & (42.77) & (983.86) & (1102.52) & (1102.52) & & (28.50) \\
    	300 & \textit{199.88} & \textbf{217.60} & 400.63 & 445.11 & 222.50 & 16845.87 & 120071.46 & 123161.75 & 123161.75 & 200.02 & 217.75 \\
    	& & (27.93) & & (43.70) & & (2167.53) & (1191.24) & (3757.26) & (3757.26) & & (27.92) \\
    	400 & \textit{199.92} & \textbf{196.61} & 627.97 & 618.38 & 222.51 & 43325.17 & 159926.82 & 161662.45 & 161662.45 & 200.00 & 196.70 \\
    	& & (27.35) & & (59.10) & & (5447.99) & (1418.60) & (3304.65) & (3304.65) & & (27.35) \\
    	500 & \textit{199.94} & \textbf{193.02} & 794.03 & 823.27 & 225.01 & 6497.32 & 199982.59 & 200043.69 & 200043.69 & 200.00 & 193.33 \\
    	& & (27.16) & & (62.69) & & (824.72) & (1550.64) & (1647.47) & (1647.47) & & (27.18) \\
    	1000 & \textit{199.97} & \textbf{185.98} & 2116.68 & 2108.78 & 224.77 & 3145.06 & 399770.90 & 400168.82 & 400168.82 & 200.03 & 186.02 \\
    	& & (27.17) & & (101.50) & & (411.60) & (2359.10) & (2934.25) & (2934.25) & & (27.17) \\
    	\bottomrule
    	\end{tabular}
 \end{table}

 \begin{table}[!h]
    	\centering
    	\caption{Sparse-robust $\alpha$ (five large coefficients equal to $10$ and other coefficients equal to $0.5$ or $-0.5$ randomly, giving $\sum_{i=1}^{n}{\alpha_i^2} = 523.75$); $X$ with all singular values equal to $1$.}
    	\label{tab:s7}
\footnotesize
	\begin{tabular}{c|cc|cc|cc|c|c|c|cc}
    	\toprule
		&\multicolumn{2}{c|}{RR} &\multicolumn{2}{c|}{PCR} & \multicolumn{2}{c|}{LASSO} & A\_LASSO & MCP & SCAD & \multicolumn{2}{c}{HS}\\
    	$p$ & SURE & SSE & SURE & SSE & SURE & SSE & SSE & SSE & SSE & SURE & SSE \\
    	\toprule
    	100 & 183.50 & 179.99 & 291.45 & 275.49 & 139.29 & 139.39 & 129.30 & 126.70 & 126.19 & \textit{131.81} & \textbf{122.60} \\
    	& & (25.31) & & (32.44) & & (20.36) & (19.20) & (18.79) & (18.83) & & (18.72) \\
    	200 & 184.47 & 196.14 & 261.65 & 277.35 & 135.93 & 150.17 & \textbf{128.76} & 129.25 & 129.58 & \textit{129.15} & 131.16 \\
    	& & (28.79) & & (33.30) & & (21.76) & (17.75) & (17.61) & (17.90) & & (18.63) \\
    	300 & 182.50 & 192.24 & 267.96 & 269.04 & 126.35 & 146.72 & 132.26 & 132.05 & 132.17 & \textit{119.09} & \textbf{128.72} \\
    	& & (25.37) & & (30.40) & & (18.96) & (17.85) & (17.87) & (17.73) & & (17.34) \\
    	400 & 184.03 & 178.58 & 311.01 & 287.95 & 145.28 & 139.68 & 128.94 & 127.57 & 127.41 & \textit{130.13} & \textbf{123.83} \\
    	& & (25.20) & & (32.98) & & (19.20) & (17.40) & (17.42) & (17.43) & & (17.02) \\
    	500 & 183.81 & 173.44 & 278.35 & 268.85 & 147.54 & 139.74 & 126.65 & 127.19 & 126.30 & \textit{132.70} & \textbf{120.78} \\
    	& & (24.08) & & (30.78) & & (19.98) & (18.36) & (18.31) & (18.17) & & (17.29) \\
    	1000 & 185.36 & 166.39 & 280.59 & 262.54 & 124.52 & 130.61 & 128.83 & 129.78 & 129.47 & \textit{119.24} & \textbf{125.49} \\
    	& & (23.16) & & (30.01) & & (18.02) & (17.50) & (17.70) & (17.62) & & (17.37) \\
    	\bottomrule
    	\end{tabular}
 \end{table}

 \begin{table}[!h]
    	\centering
    	\caption{Null $\alpha$ ($\sum_{i=1}^{n} \alpha_i^2 = 0$); $X$ with all singular values equal to $1$.}
    	\label{tab:s8}
\footnotesize
	\begin{tabular}{c|cc|cc|cc|c|c|c|cc}
    	\toprule
		&\multicolumn{2}{c|}{RR} &\multicolumn{2}{c|}{PCR} & \multicolumn{2}{c|}{LASSO} & A\_LASSO & MCP & SCAD & \multicolumn{2}{c}{HS}\\
    	$p$ & SURE & SSE & SURE & SSE & SURE & SSE & SSE & SSE & SSE & SURE & SSE \\
    	\toprule
    	100 & 94.70 & \textbf{100.13} & 97.63 & 102.62 & \textit{94.54} & 100.15 & \textbf{100.13} & 100.15 & 100.15 & 99.92 & 101.44 \\
    	& & (14.71) & & (15.37) & & (14.69) & (14.70) & (14.69) & (14.69) & & (15.06) \\
    	200 & 115.52 & 103.43 & 111.09 & 118.11 & \textit{109.16} & 122.81 & 100.49 & 112.22 & \textbf{100.80} & 116.62 & 106.72 \\
    	& & (14.80) & & (16.91) & & (17.79) & (14.37) & (16.20) & (14.53) & & (15.47) \\
    	300 & 98.74 & 100.49 & 99.45 & 113.35 & \textit{96.40} & 103.03 & \textbf{100.46} & 103.03 & 103.03 & 103.10 & 102.40 \\
    	& & (14.80) & & (16.49) & & (15.04) & (14.82) & (15.04) & (15.04) & & (14.87) \\
    	400 & 96.78 & \textbf{97.99} & 103.88 & 102.24 & \textit{94.02} & 103.08 & \textbf{97.99} & 101.71 & 103.17 & 99.64 & 101.01 \\
    	& & (14.49) & & (15.12) & & (14.96) & (14.50) & (14.81) & (14.97) & & (14.78) \\
    	500 & 88.55 & \textbf{99.97} & 89.06 & 100.91 & \textit{87.74} & 100.65 & 99.98 & 100.65 & 100.65 & 93.63 & 101.53 \\
    	& & (14.81) & & (14.72) & & (14.87) & (14.83) & (14.87) & (14.87) & & (14.98) \\
    	1000 & 88.87 & \textbf{100.94} & 91.96 & 107.30 & \textit{88.45} & 101.14 & 100.95 & 101.62 & 101.14 & 94.34 & 102.26 \\
    	& & (14.17) & & (15.40) & & (14.14) & (14.17) & (14.30) & (14.14) & & (14.48) \\
    	\bottomrule
    	\end{tabular}
 \end{table}

\begin{table}[!h]
    	\centering
    	\caption{Dense $\alpha$ (all coefficients equal to 2, giving $\sum_{i=1}^{n}{\alpha_i^2} = 400$); $X$ with all singular values equal to $1$.}
    	\label{tab:s9}
\footnotesize
	\begin{tabular}{c|cc|cc|cc|c|c|c|cc}
    	\toprule
		&\multicolumn{2}{c|}{RR} &\multicolumn{2}{c|}{PCR} & \multicolumn{2}{c|}{LASSO} & A\_LASSO & MCP & SCAD & \multicolumn{2}{c}{HS}\\
    	$p$ & SURE & SSE & SURE & SSE & SURE & SSE & SSE & SSE & SSE & SURE & SSE \\
    	\toprule
    	100 & \textit{177.89} & \textbf{183.69} & 220.87 & 200.41 & 203.16 & 307.44 & 502.55 & 505.43 & 505.43 & 200.80 & 204.16 \\
    	& & (25.16) & & (26.35) & & (43.19) & (41.53) & (43.20) & (43.20) & & (27.46) \\
    	200 & \textit{181.88} & \textbf{188.39} & 207.49 & 239.95 & 214.11 & 255.88 & 499.88 & 498.84 & 498.84 & 205.16 & 217.46 \\
    	& & (27.10) & & (33.45) & & (35.61) & (41.69) & (42.90) & (42.90) & & (30.81) \\
    	300 & \textit{176.64} & \textbf{193.53} & 215.03 & 205.00 & 196.19 & 250.76 & 496.84 & 497.60 & 495.93 & 199.33 & 212.19 \\
    	& & (25.76) & & (26.70) & & (32.66) & (45.74) & (47.41) & (45.90) & & (27.68) \\
    	400 & \textit{174.62} & \textbf{195.83} & 248.90 & 249.51 & 209.80 & 221.41 & 495.21 & 494.17 & 494.17 & 198.84 & 206.49 \\
    	& & (26.36) & & (32.46) & & (29.85) & (40.50) & (40.89) & (40.89) & & (28.41) \\
    	500 & \textit{179.13} & \textbf{173.13} & 202.78 & 192.97 & 214.36 & 201.27 & 501.67 & 503.19 & 503.19 & 205.16 & 193.40 \\
    	& & (23.18) & & (25.08) & & (25.63) & (38.90) & (38.93) & (38.93) & & (25.60) \\
    	1000 & \textit{179.32} & \textbf{173.71} & 225.50 & 195.46 & 209.35 & 248.20 & 503.44 & 507.29 & 507.29 & 204.51 & 194.73 \\
    	& & (24.14) & & (25.78) & & (30.87) & (41.48) & (41.38) & (41.38) & & (26.67) \\
    	\bottomrule
    	\end{tabular}
 \end{table}

\end{document}